\documentclass{amsart}
\usepackage[cp850]{inputenc}

\newtheorem{theorem}{Theorem}[section]
\newtheorem{claim}{Claim}[section]
\newtheorem{lemma}[theorem]{Lemma}
\newtheorem{proposition}[theorem]{Proposition}
\newtheorem{corollary}[theorem]{Corollary}

\theoremstyle{definition}
\newtheorem{definition}[theorem]{Definition}

\theoremstyle{remark}

\numberwithin{equation}{section}

\newcommand{\g}[2]{\mbox{$\langle #1 ,#2 \rangle$}}
\newcommand{\mbar}{\mbox{$\overline{M}$}}
\newcommand{\m}{\mbox{$M$}}
\newcommand{\nablabar}{\mbox{$\overline{\nabla}$}}

\newcommand{\xm}{\mbox{$\mathcal{X}(\m)$}}
\newcommand{\x}{\mbox{$\psi:\m^n\rightarrow\mbar^{n+1}$}}
\newcommand{\rf}[1]{\mbox{(\ref{#1})}}
\newcommand{\rl}[1]{{~\ref{#1}}}

\def\Div{\mathrm{div}_{M}\ }
\def\Divbar{\mathrm{div}_{\overline{M}}\ }
\def\fle{\rightarrow}
\def\beq{\begin{equation}}
\def\eeq{\end{equation}}
\def\Rn{\mbox{$\mathbb R^{n+1}$}}
\def\Hn{\mbox{$\mathbb H^{n+1}$}}
\def\Sn{\mbox{$\mathbb S^{n+1}$}}
\def\Snp{\mbox{$\mathbb S^{n+1}_{+}$}}

\newcommand{\xR}{\mbox{$\psi:\m^n\rightarrow\Rn$}}
\newcommand{\xS}{\mbox{$\psi:\m^n\rightarrow\Sn$}}
\newcommand{\xH}{\mbox{$\psi:\m^n\rightarrow\Hn$}}

\begin{document}

\title[Constant higher order mean curvature hypersurfaces]
{Constant higher order mean curvature hypersurfaces in Riemannian spaces}

\author{Luis J. Al\'\i as}
\address{Departamento de Matem\'aticas, Universidad de Murcia, E-30100 Espinardo, Murcia,
Spain}
\email{ljalias@um.es}
\thanks{L.J. Al\'\i as was partially supported by MCyT and Fundaci\'on S\'eneca, Spain}

\author{Jorge H.S. de Lira}
\address{Departamento de Matem\'atica, Universidade Federal do Ceará,  Campus do Pici, 
60455-760 Fortaleza-Ce, Brazil}
\email{jherbert@mat.ufc.br}
\thanks{J.H.S. de Lira  was partially supported by CNPq, Brazil}

\author{J. Miguel Malacarne}
\address{Departamento de Matem\'atica, Universidade Federal do Esp\'\i rito Santo, 29060-900
Vit\'oria-ES, Brazil}
\email{jmiguel@cce.ufes.br}
\thanks{J.M. Malacarne was partially supported by PICD/CAPES, Brazil}
\thanks{Some of the main results contained in this paper were previously announced in the 
International Conference on Differential Geometry which was held to 
honour the 60th birthday of Professor A.M. Naveira in 
Valencia, Spain, in July 2001, as well as in the XII Escola de Geometria 
Diferencial held in Goiânia, Goiás, Brazil, in July 2002}


\date{\bf July 16, 2002}



\begin{abstract}
It is still an open question whether a compact embedded hypersurface in the Euclidean space \Rn\ 
with constant mean curvature and spherical boundary is necessarily a hyperplanar ball or a 
spherical cap, even in the simplest case of surfaces in $\mathbb{R}^3$. In a recent paper 
\cite{AM} the first and third authors have shown that this is true for the case of 
hypersurfaces  in \Rn\ with constant scalar curvature, and more generally, hypersurfaces with 
constant higher order $r$-mean curvature, when $r\geq 2$. In this paper we deal 
with some aspects of the classical problem above, by considering it in a 
more general context. Specifically, our starting general ambient space is an 
orientable Riemannian manifold \mbar, where we will consider a general geometric configuration 
consisting of an immersed hypersurface into \mbar\ with boundary on an oriented hypersurface $P$ 
of \mbar. For such a geometric configuration, we study the relationship between the geometry of 
the hypersurface along its boundary and the geometry of its boundary as a hypersurface of $P$, 
as well as the geometry of $P$ as a hypersurface of \mbar. 
Our approach allows us to derive, among others, interesting results for the case where the 
ambient space has constant curvature (the Euclidean space \Rn, the hyperbolic space \Hn, and 
the sphere \Sn). In particular, we are able to extend the symmetry results given in \cite{AM} to 
the case of hypersurfaces with constant higher order $r$-mean curvature in the hyperbolic space 
and in the sphere.
\end{abstract}

\maketitle

\section{Introduction}
An old problem in classical differential geometry consists on finding all 
compact surfaces in Euclidean space $\mathbb{R}^3$ with constant mean 
curvature and circular boundary. As is well known, a circle $C$ in 
$\mathbb{R}^3$ is the boundary of two spherical caps with constant mean curvature $H$ for any 
positive number $H$, less than or equal to the inverse of the radius of the circle $C$. 
A natural question to ask \cite{B} is whether a compact constant mean 
curvature surface in $\mathbb{R}^3$ which is bounded by a circle is necessarily a spherical cap 
or a flat disc. Actually, a constant mean curvature surface with circular boundary is the 
mathematical model of a soap bubble which has its boundary on a round hoop, and the surfaces we 
almost always observe are spherical caps, so that it is natural to ask if these are the only 
solutions. In \cite{Ka} Kapouleas gave a negative answer to this question by showing that 
there exist examples of higher genus compact, non-spherical immersed surfaces with
constant mean curvature in $\mathbb{R}^3$ bounded by a circle. 
However, it has been conjectured that there must be a positive answer to this question 
if one requires in addition that the surface has genus zero or that it is embedded 
\cite{BEMR}. 

In recent years, several authors have obtained some partial answers to 
these problems. For instance, Barbosa \cite{B1,B2} proved that the only 
compact immersed surfaces with constant mean curvature $H\neq 0$ and circular boundary which are 
contained either in a sphere or in a cilinder of radius $1/|H|$ are the spherical caps. 
On the other hand, in the genus zero case the first author, jointly with López and Palmer, 
has showed that the only stable constant mean curvature immersed surfaces of disc type which are 
bounded by a circle are spherical caps \cite{ALP} (see also \cite{BJ} for another 
characterization of spherical caps as the only stable examples, given by Barbosa and Jorge under 
a stronger idea of stability). 

It is clear that this classical question can be stated in a more general 
context as follows. Let $\Sigma^{n-1}$ be a compact $(n-1)$-dimensional 
submanifold contained in a hyperplane $\Pi\subset\Rn$, and let $\m^n$ be an $n$-dimensional
connected orientable manifold with smooth boundary $\partial\m$. As usual, \m\ is said to be a
hypersurface of \Rn\ with boundary $\Sigma$ if there exists an immersion \xR\ such that the 
immersion $\psi$ restricted to the boundary $\partial\m$ is a diffeomorphism 
onto $\Sigma$. In this context, the classical question above consists on 
finding the compact hypersurfaces in \Rn\ with constant mean curvature 
whose boundary $\Sigma$ is a round $(n-1)$-sphere. At this point, it is interesting to 
recall that a classical result by Alexandrov \cite{Al} states that round 
spheres are the only \textit{closed} hypersurfaces with constant mean curvature 
which are embedded in Euclidean space \Rn\ (here by closed we mean compact 
and \textit{without} boundary). More recently, Alexandrov theorem was 
extended by Ros to the case of constant scalar curvature \cite{Ro1}, and more generally to the 
case of hypersurfaces with constant higher order mean curvature \cite{Ro2}, showing that round 
spheres are the only closed embedded hypersurfaces with constant $r$-mean curvature in 
\Rn.

As for the case of non-empty boundary, in \cite{Ko} Koiso gave a new interpretation of the 
problem by studying 
under what conditions the symmetries of the boundary $\Sigma\subset\Pi$ of a non-zero constant 
mean curvature hypersurface \m\ in \Rn\ are inherited by the whole hypersurface. She 
showed that this necessarily occurs when the hypersurface \m\ is embedded and 
it does not intersect the outside of $\Sigma$ in $\Pi$; as a consequence, if the boundary 
$\Sigma$ is a round $(n-1)$-sphere, then \m\ is symmetric with respect to every 
hyperplane through the center of $\Sigma$ which is orthogonal to $\Pi$, 
and hence \m\ must be a spherical cap. 
Related to Koiso's symmetry theorem, Brito, S\'a Earp, Meeks and Rosenberg \cite{BEMR}
also showed that when $\Sigma$ is strictly convex and \m\ is embedded and transverse to $\Pi$ 
along the boundary $\partial\m$, then \m\ is entirely contained in one of the half-spaces of 
\Rn\ determined by $\Pi$ and, therefore, the so called \textit{Alexandrov reflection technique} 
\cite{Al} implies that \m\ inherites all the symmetries of $\Sigma$. 
In particular, if $\Sigma$ is a round sphere, then \m\ must be a spherical cap.
Here, transversality means that the hypersurface \m\ is never tangent to the hyperplane $\Pi$ 
along its boundary. In what follows, we will use the term \textit{symmetry 
result} to refer to a result of this type.

The technique introduced in \cite{BEMR} makes an extensive use of two 
essential ingredients, the abovementioned Alexandrov reflection technique,
and an integral formula first found by Kusner \cite{Ku}, which is now known as 
the \textit{flux formula}. This fact indicates that the symmetry result in 
\cite{BEMR} can be extended from two new viewpoints: by considering 
constant mean curvature hypersurfaces in other space forms; or by 
considering the case of hypersurfaces with constant higher order $r$-mean 
curvature. From the first point of view, Nelli and Rosenberg \cite{NR} 
studied the case of hypersurfaces in hyperbolic space \Hn, and, more recently, Lira 
\cite{Lira} considered the case of hypersurfaces in the sphere \Sn, 
establishing corresponding symmetry results for the case of constant mean curvature. On the 
other hand, in \cite{Ro} Rosenberg established a version of the flux formula for hypersurfaces 
with constant higher order $r$-mean curvature in Euclidean space \Rn, and applied it to extend 
the symmetry result given in \cite{BEMR} to the case of the higher order $r$-mean curvatures.

In this paper, we will deal with some aspects of the classical problem above. Our initial 
strategy is to study this problem in a more general context. Specifically, our general ambient 
space will be an $(n+1)$-dimensional connected orientable Riemannian manifold \mbar, where we 
will consider the following geometric configuration (for the details, see 
Section\rl{configuration}). Let us fix $P^n\subset\mbar$ an orientable connected totally geodesic 
hypersurface in \mbar, and let $\Sigma^{n-1}\subset P$ be an orientable $(n-1)$-dimensional 
compact embedded submanifold contained in $P^n$. Consider $\m^n$ an $n$-dimensional connected 
orientable manifold with smooth boundary $\partial\m$. Then, \m\ is said to be a
hypersurface of \mbar\ with boundary $\Sigma$ if there exists an immersion \x\ such that the 
immersion $\psi$ restricted to the boundary $\partial\m$ is a diffeomorphism onto $\Sigma$.
From this geometric configuration, the following question, closely related to the symmetry 
problem, naturally arises: 

\begin{it}
How is the geometry of \m\ along its boundary $\partial\m$ related to 
the geometry of the inclusion $\Sigma\subset P$?
\end{it}

A first partial answer to this question is given by the following 
expression, which holds along the boundary $\partial\m$ and for every 
$1\leq r\leq n-1$ (see Corollary\rl{c1}),
\[
\g{T_r\nu}{\nu}=(-1)^rs_r\g{\xi}{\nu}^{r}.
\]
Here $T_r$ stands for the $r$-th classical Newton transformation associated to 
the second fundamental form on \m\ (see Section\rl{Newton} for the 
details), $\nu$ is the outward pointing unit conormal vector field along $\partial\m$, $\xi$ is 
the unitary normal field of $P\subset\mbar$, and $s_r=s_r(\tau_1,\ldots,\tau_{n-1})$
is the $r$-th elementary symmetric function of $\tau_1,\ldots,\tau_{n-1}$, the principal
curvatures of $\Sigma\subset P$ with respect to the outward pointing unitary
normal. As a first consequence of this expression, we obtain a very strong 
relationship between the transversality of \m\ with respect to $P$ along 
the boundary $\partial\m$, and the ellipticity on \m\ of the $r$-th Newtom 
tranformation $T_r$, that is, the positivity of the quadratic form 
associated to $T_r$. This fact, along with Theorem 7.3 in \cite{Ro}, 
allows us to state the following symmetry theorem for hypersurfaces in 
\Rn\ (Theorem\rl{calotaRn}):

\begin{it}
Let $\Sigma$ be an strictly convex compact $(n-1)$-dimensional submanifold in a
hyperplane $\Pi\subset\Rn$, and let \xR\ be a compact embedded hypersurface with boundary
$\Sigma$. Let us assume that for a given $2\leq r\leq n$, the $r$-mean
curvature $H_r$ of \m\ is a nonzero constant . Then \m\ has all the
symmetries of $\Sigma$. In particular, if the boundary $\Sigma$ is a round
$(n-1)$-sphere of \Rn, then \m\ is a spherical cap.
\end{it}

As a consequence, we can conclude that the conjecture of the spherical cap 
\cite{BEMR} is true for the case of embedded hypersurfaces with constant 
$r$-mean curvature in \Rn, when $r\geq 2$. This includes, in particular, 
the case of constant scalar curvature, when $r=2$ \cite{AM}.

In order to extend this symmetry result to the case of hypersurfaces in 
hyperbolic space and hypersurfaces in the sphere, it is necessary to 
establish a certain \textit{flux formula}, which is one of the key 
ingredients of the used techniques. For that reason, 
Section\rl{sectionflux} is devoted to derive a general flux formula for 
the considered geometric configuration in the case where the Riemannian 
ambient space \mbar\ is equipped with a conformal vector field (Proposition\rl{propflux}). 
Our general flux formula becomes specially simple when the ambient space 
has constant sectional curvature, and the conformal vector field is indeed a 
Killing vector field. In that case, we are able to extend the flux formula 
given by Rosenberg in \cite[Theorem 7.2]{Ro} to the case of the other 
space forms, as follows (Corollary\rl{coroflux}):

\begin{it}
Let \x\ be an immersed compact orientable hypersurface with boundary $\partial\m$, and
let $D^n$ be a compact orientable hypersurface with boundary $\partial D=\partial\m$. Assume 
that $\m\cup D$ is an oriented $n$-cycle of \mbar,
and let ${\bf N}$ and $n_D$ be the unit normal fields which orient \m\ and $D$, respectively.
Assume that \mbar\ has constant sectional curvature.
If the $r$-mean curvature $H_r$ is constant, $1\leq r\leq n$, then for every Killing vector
field $Y\in\mathcal{X}(\mbar)$ the following flux formula holds
\[
\oint_{\partial M}\g{T_{r-1}\nu}{Y}ds=-r\binom{n}{r}H_r\int_D\g{Y}{n_D}dD,
\]
where $\nu$ is the outward pointing conormal to \m\ along $\partial\m$.
\end{it}

As first applications of our general flux formula, we derive some interesting 
estimates for the volume of minimal hypersurfaces
with boundary on a geodesic sphere 
of the ambient space, in the case where the ambient space is the Euclidean 
space (Corollary\rl{coroflux1}), the hyperbolic space 
(Corollary\rl{coroflux2}), or the sphere (Corollary\rl{coroflux3}). 

On the other hand, and as another application of our flux formula and the expression for 
$\g{T_r\nu}{\nu}$ given in Corollary\rl{c1} (see above), we establish in Section\rl{estimating} 
some interesting estimates for the constant $r$-mean curvature in terms of 
the geometry of the boundary. Specifically, when the ambient space is the 
Euclidean space we obtain the following (Theorem\rl{estimate1}):

\begin{it}
Let $\Sigma$ be an orientable $(n-1)$-dimensional compact submanifold in a hyperplane
$P\subset\Rn$, and let \xR\ be an orientable immersed compact (connected)
hypersurface with boundary $\Sigma=\psi(\partial\m)$ and constant $r$-mean curvature $H_r$,
$1\leq r\leq n$. Then
\[
0\leq |H_r|\leq\frac{1}{n\ \mathrm{vol}(D)}\oint_{\partial M}|h_{r-1}|ds,
\]
where $h_{r-1}$ stands for the $(r-1)$-mean curvature of $\Sigma\subset
P$, and $D$ is the domain in $P$ bounded by $\Sigma$. In particular, when
$\Sigma$ is a round $(n-1)$-sphere of radius $\varrho$ it follows that
\[
0\leq|H_r|\leq\frac{1}{\varrho^r}.
\]
\end{it}

This estimate is the natural generalization of an estimate first obtained by Barbosa 
in the case of constant mean curvature ($r=1$) \cite{B1}. On the other 
hand, when the ambient space is the hyperbolic space, our estimate reads 
as follows (Theorem\rl{estimate2}):

\begin{it}
Let $\Sigma$ be an orientable $(n-1)$-dimensional compact submanifold contained in a totally 
geodesic hyperplane $P\subset\Hn$, and let \xH\ be an orientable immersed compact connected
hypersurface with boundary $\Sigma=\psi(\partial\m)$ and constant $r$-mean curvature $H_r$,
$1\leq r\leq n$. Then
\[
0\leq |H_r|\leq\frac{C}{n\ \mathrm{vol}(D)}\oint_{\partial M}|h_{r-1}|ds.
\]
Here $h_{r-1}$ stands for the $(r-1)$-mean curvature of $\Sigma\subset P$, $D$ is the domain in 
$P$ bounded by $\Sigma$, and $C=\max_{\Sigma}\cosh{\tilde{\varrho}}\geq 1$, where 
$\tilde{\varrho}(p)$ is the geodesic distance along $P$ between a fixed arbitrary point 
$a\in\mathrm{int}(D)$ and $p$. 
In particular, when $\Sigma$ is a geodesic sphere in $P$ of geodesic radius $\varrho$, it 
follows that
\[
0\leq|H_r|\leq\coth^{r}\varrho.
\]
\end{it}

Similarly, for the case of hypersurfaces in the sphere, our estimate states 
the following (Theorem\rl{estimate3}):

\begin{it}
Let $\Sigma$ be an orientable $(n-1)$-dimensional compact submanifold contained in an open 
totally geodesic hemisphere $P_+\subset\Sn$, and let \xS\ be an orientable immersed compact 
connected hypersurface with boundary $\Sigma=\psi(\partial\m)$ and constant $r$-mean curvature 
$H_r$, $1\leq r\leq n$. Then
\[
0\leq |H_r|\leq\frac{C}{n\ \mathrm{vol}(D)}\oint_{\partial M}|h_{r-1}|ds.
\]
Here $h_{r-1}$ stands for the $(r-1)$-mean curvature of $\Sigma\subset P$, $D$ is the domain in 
$P_+$ bounded by $\Sigma$, and $C=\max_{\Sigma}\cos{\tilde{\varrho}}/\min_D\cos{\tilde{\varrho}}$, 
where $\tilde{\varrho}(p)$ is the geodesic distance along $P_+$ between a fixed arbitrary point 
$a\in\mathrm{int}(D)$ and $p$. 
In particular, when $\Sigma$ is a geodesic sphere in $P_+$ of geodesic radius $\varrho<\pi/2$, 
it follows that
\[
0\leq|H_r|\leq\cot^{r}\varrho.
\]
\end{it}

Finally, the two remaining sections of the paper are devoted to the 
extension of our symmetry results to the case of hypersurfaces in the hyperbolic space and 
hypersurfaces in the sphere. Specifically, in Section\rl{hyperbolic} we 
obtain the following symmetry result for hypersurfaces in hyperbolic space 
(Theorem\rl{simetriaHn}):

\begin{it}
Let $\Sigma^{n-1}$ be an strictly convex compact $(n-1)$-dimensional 
(connected) submanifold of a totally geodesic hyperplane $P^n\subset\Hn$, and let 
$M^n\subset\Hn$ be a compact (connected) embedded  hypersurface with boundary $\Sigma$. 
Let us assume that for a given $2\leq r\leq n$, the $r$-mean curvature 
$H_r$ of $M$ is a nonzero constant. Then $M$ has all the symmetries of $\Sigma$. 
In particular, when the boundary $\Sigma$ is a geodesic sphere in $P^n\subset\Hn$, then  
$M$ is a spherical cap.
\end{it}

As a consequence,  we can conclude, as in the Euclidean case, that the conjecture of the 
spherical cap is true for the case of embedded hypersurfaces with constant $r$-mean curvature in 
hyperbolic space, when $r\geq 2$. Finally, in the case of hypersurfaces in 
the sphere, we state the following symmetry result (Theorem\rl{simetriaSn}):

\begin{it}
Let $\Sigma^{n-1}$ be a convex $(n-1)$-dimensional submanifold of a totally geodesic $n$-sphere 
$P^n\subset\Sn$, and let $M^n\subset\Sn$ be a compact (connected) embedded hypersurface with 
boundary $\Sigma$. Let us assume that $M$ is contained in an open hemisphere \Snp, and that the 
$r$-mean curvature $H_r$ 
of $M$ is a nonzero constant, for a given $2\leq r\leq n$.  Then $M$ has all the symmetries of 
$\Sigma$. In particular, when the boundary $\Sigma$ is a geodesic sphere in $P^n\subset\Sn$, 
then $M$ is a spherical cap.
\end{it}

In particular, the only compact embedded hypersurfaces in \Snp\ with constant $r$-mean curvature $H_r\neq 0$
(with $2\leq r\leq n$) and spherical boundary are the spherical caps.

\section{Preliminaries}
\label{preliminaries}
Throughout this paper, $\mbar^{n+1}$ will denote an $(n+1)$-dimensional connected orientable Riemannian
manifold, and $\langle \;,\;\rangle$ and \nablabar\ will stand for its Riemannian metric and
its Levi-Civita connection, respectively. Let $\m^n$ be an $n$-dimensional
connected orientable manifold with smooth boundary $\partial\m$; \m\ is said to be a
hypersurface of \mbar\ if there exists an isometric immersion \x. In that case, since \m\ and
\mbar\ are both orientable, we may choose along $\psi(M)$ a globally defined unit normal vector field
${\bf N}$, and we may assume that \m\ is oriented by ${\bf N}$. If $\nabla$ denotes the Levi-Civita
connection on \m\, then the Gauss and Weingarten formulae for the immersion are given,
respectively by
\beq
\label{G}
\nablabar_VW=\nabla_VW+\g{AV}{W}{\bf N},
\eeq
and
\beq
\label{W}
A(V)=-\nablabar_V{\bf N},
\eeq
for all tangent vector fields $V,W\in\xm$.

Here $A:\xm\fle\xm$ defines the shape operator (or the second fundamental form) of the
hypersurface with respect to ${\bf N}$. The curvature tensor $R$ of the hypersurface \m\ is described
in terms of $A$ and the curvature tensor $\overline{R}$ of the ambient space \mbar\ by the
so called Gauss equation, which can be written as
\beq
\label{gausseq}
R(U,V)W=(\overline{R}(U,V)W)^\top+\g{AU}{W}AV-\g{AV}{W}AU
\eeq
for all tangent vector fields $U,V,W\in\xm$, where ${}^\top$ denotes
projection on \xm. Observe that our criterion here for the definition of the curvature tensor is
the one in \cite{ONe}. On the other hand, the Codazzi equation of the hypersurface describes the
normal component of $\overline{R}(U,V)W$ in terms of the derivative of the shape operator, and
it is given by
\beq
\label{codazzi}
\g{\overline{R}(U,V)W}{{\bf N}}=\g{(\nabla_VA)U-(\nabla_UA)V}{W}
\eeq
where $\nabla_UA$ denotes the covariant derivative of $A$. In particular, when the ambient
space has constant sectional curvature, then $\overline{R}(U,V)W$ is tangent to \m\ for every
$U,V,W\in\xm$, and \rf{codazzi} becomes
\beq
\label{codazzicc}
(\nabla_VA)U=(\nabla_UA)V.
\eeq

As is well known, $A$ is a self-adjoint linear operator in each tangent plane $T_pM$, and
its eigenvalues $\kappa_1(p), \ldots, \kappa_n(p)$ are the principal curvatures
of the hypersurface. Associated to the shape operator there are $n$
algebraic invariants given by
\[
S_r(p)=\sigma_r(\kappa_1(p), \ldots, \kappa_n(p)), \quad 1\leq r\leq n.
\]
where $\sigma_r:\mathbb R^n\to\mathbb R$ is the elementary symmetric functions in $\mathbb R^n$ given by 
\[
\sigma_r(x_1,\ldots, x_n)=\sum_{i_1<\cdots<i_r}x_{i_1}\ldots x_{i_n}.
\]
Observe that the characteristic polynomial of $A$ can be writen in terms
of the $S_r$'s as
\beq
\label{poly}
\det(tI-A)=\sum_{^r=0}^n(-1)^rS_rt^{n-r}.
\eeq
The $r$-mean curvature $H_r$ of the hypersurface is then defined by
\[
\binom{n}{r}H_r=S_r.
\]
In particular, when $r=1$ then $H_1=(1/n)\mathrm{trace}(A)=H$ is the mean curvature of \m,
which is the main extrinsic curvature of the hypersurface. On the other hand, when $r=2$, $H_2$
defines a geometric quantity which is related to the (intrinsic) scalar curvature of the
hypersurface. Indeed, it follows from the Gauss equation \rf{gausseq} that the Ricci curvature
of \m\ is given by
\[
\mathrm{Ric}(U,V)=\mathrm{\overline{Ric}}(U,V)-\g{\overline{R}(U,{\bf N})V}{{\bf N}}
+nH\g{AU}{V}-\g{AU}{AV},
\]
for $U,V\in\xm$, where $\mathrm{\overline{Ric}}$ stands for the Ricci curvature
of the ambient space \mbar. Therefore, the scalar curvature $S$ of the hypersurface \m\
is
\[
S=\mathrm{trace}(\mathrm{Ric})=\overline{S}-2\ \overline{\mathrm{Ric}}({\bf N},{\bf N})+n(n-1)H_2.
\]
For instance, if the ambient space has constant sectional curvature $\overline{c}$ then
\beq
\label{scalar}
S=n(n-1)(\overline{c}+H_2).
\eeq

\section{The Newton transformations}
\label{Newton}
The classical Newton transformations $T_r:\xm\fle\xm$ are defined inductively from $A$ by
\[
T_0=I \quad \mathrm{and} \quad T_r=S_rI-AT_{r-1}, \quad 1\leq r\leq n,
\]
where $I$ denotes the identity in \xm, or equivalently by
\[
T_r=S_rI-S_{r-1}A+\cdots+(-1)^{r-1}S_1A^{r-1}+(-1)^{r}A^{r}.
\]
Note that by the Cayley-Hamilton theorem, we have $T_n=0$.

Let us recall that each $T_r$ is also a self-adjoint linear operator in each tangent plane $T_pM$ which
commutes with $A$. Indeed, $A$ and $T_r$ can be simultaneously
diagonalized; if $\{ e_1, \ldots, e_n\}$ are the eigenvectors of $A$
corresponding to the eigenvalues $\kappa_1(p), \ldots, \kappa_n(p)$, respectively, then they
are also the eigenvectors of $T_r$ corresponding to the eigenvalues of
$T_r$, and $T_r(e_i)=\mu_{i,r}(p)e_i$ with
\[
\mu_{i,r}(p)=\frac{\partial \sigma_{r+1}}{\partial x_i}(\kappa_1(p), \ldots, \kappa_n(p))
=\sum_{i_1<\cdots<i_r,i_j\neq i}\kappa_{i_1}(p)\cdots\kappa_{i_r}(p),
\]
for every $1\leq i\leq n$. From here it can be easily seen that
\begin{eqnarray}
\mathrm{trace}(T_r) & = & (n-r)S_r=c_rH_r, \label{traces} \\
\mathrm{trace}(AT_r) & = & (r+1)S_{r+1}=c_rH_{r+1}, \label{tracesbis}
\end{eqnarray}
where $c_r=(n-r)\binom{n}{r}=(r+1)\binom{n}{r+1}$.
For the details, we refer the reader to the classical paper by Reilly
\cite{Re} (see also \cite{Ro} for a more accesible modern treatment by Rosenberg).

On the other hand, the divergence of $T_r$ is defined by
\[
\Div T_r=\mathrm{trace}(\nabla T_r)=
\sum_{i=1}^n(\nabla_{e_i}T_r)(e_i),
\]
where $\{ e_1, \ldots, e_n \}$ is a local orthonormal frame on \m. Below
we will compute $\Div T_r$, which will be necessary for its later use.
\begin{lemma}
\label{divlema}
The divergence of the Newton transformations $T_r$ are given by the following inductive formula,
\begin{eqnarray}
\label{div}
\Div T_0 & = & 0 \\
\nonumber \Div T_r & = &
-A(\Div T_{r-1})-\sum_{i=1}^{n}(\overline{R}({\bf N},T_{r-1}e_i)e_i)^\top,
\end{eqnarray}
where $\overline{R}$ stands for the curvature tensor of \mbar, and
$(\overline{R}({\bf N},V)W)^\top$ denotes the tangential component of
$\overline{R}({\bf N},V)W$. Equivalently, for every tangent field $V\in\xm$ it follows
\begin{equation}
\label{otradiv}
\g{\Div T_r}{V}=\sum_{j=1}^{r}\sum_{i=1}^{n}
\g{\overline{R}({\bf N},T_{r-j}e_i)e_i}{A^{j-1}V}.
\end{equation}
\end{lemma}
The expression \rf{otradiv} has been recently obtained also by Lima
in \cite{Li}, using a very different argument to ours.
\begin{proof}
It is clear that $\Div T_0=\Div I=0$.
When $r\geq 1$, from the inductive definition of $T_r$ we have for $V,W\in\xm$
\begin{eqnarray*}
(\nabla_VT_r)W & = & \g{\nabla S_r}{V}W-\nabla_V(AT_{r-1})W\\
{} & = & \g{\nabla S_r}{V}W-(\nabla_VA)(T_{r-1}W)-A((\nabla_VT_{r-1})W),
\end{eqnarray*}
so that
\[
\Div T_r=\sum_{i=1}^n(\nabla_{e_i}T_r)(e_i)=\nabla S_r-
\sum_{i=1}^n(\nabla_{e_i}A)(T_{r-1}e_i)-A(\Div T_{r-1}).
\]
Using now the Codazzi equation \rf{codazzi} we get for $V\in\xm$
\begin{eqnarray*}
\g{(\nabla_{e_i}A)(T_{r-1}e_i)}{V} & = & \g{(\nabla_{e_i}A)V}{T_{r-1}e_i} \\
{} & = & \g{(\nabla_{V}A)e_i}{T_{r-1}e_i}+\g{\overline{R}(V,e_i)T_{r-1}e_i}{{\bf N}} \\
{} & = & \g{T_{r-1}((\nabla_{V}A)e_i)}{e_i}+\g{\overline{R}({\bf N},T_{r-1}e_i)e_i}{V}.
\end{eqnarray*}
Therefore,
\begin{eqnarray}
\label{eqdiv}
\g{\Div T_r}{V} & = & \g{\nabla S_r}{V}-\mathrm{trace}(T_{r-1}\nabla_VA)\\
\nonumber {} & - & \sum_{i=1}^{n}\g{\overline{R}({\bf N},T_{r-1}e_i)e_i}{V}-\g{A(\Div T_{r-1})}{V}.
\end{eqnarray}
Using now equation (4.4) in \cite{Ro} we have that
\[
\mathrm{trace}(T_{r-1}\nabla_VA)=\g{\nabla S_r}{V},
\]
which jointly with \rf{eqdiv} gives \rf{div}. Finally, equation \rf{otradiv}
follows easily from \rf{div} by an inductive argument.
\end{proof}
In particular, when the ambient Riemannian space \mbar\ has constant sectional curvature, then
$(\overline{R}({\bf N},V)W)^\top=0$ for every tangent vector fields $V,W\in\xm$ and
equation \rf{otradiv} implies that $\Div T_r=0$ for every $r$.
\begin{corollary}
\label{divcoro}
When the ambient Riemannian space \mbar\ has constant sectional curvature, then the
Newton transformations are divergence-free: $\Div T_r=0$ for each $r$.
\end{corollary}

\section{A geometric configuration}
\label{configuration}
Throughout this paper, we will be particularly interested in the following
geometric configuration, which is suggested by the classical question stated in the Introduction. 
Let $P^n\subset\mbar$ be an orientable connected hypersurface in
\mbar, and let $\Sigma^{n-1}\subset P$ be an orientable $(n-1)$-dimensional compact embedded
submanifold contained in $P^n$. Let \x\ be an orientable compact connected hypersurface in \mbar\
with smooth boundary $\partial\m$. As usual, \m\ is said to be a
hypersurface \textit{with boundary} $\Sigma$ if the immersion $\psi$
restricted to the boundary $\partial\m$ is a diffeomorphism onto $\Sigma$.
The following question naturally arises from this geometric configuration:
\begin{quote}
\textit{How is the geometry of \m\ along its boundary $\partial\m$ related to the
geometry of the inclusion $\Sigma\subset P$ and the inclusion
$P\subset\mbar$?}
\end{quote}

In what follows, we will study this question. Let us start by choosing the orientation of this
configuration. Let us consider the hypersurface \m\ oriented by
a globally defined unit normal vector field ${\bf N}$. The orientation of \m\
induces a natural orientation on its boundary as follows: given a point
$p\in\partial\m$, a basis $\{ v_1, \ldots, v_{n-1} \}$ for $T_p(\partial\m)$ is said to be
positively oriented if $\{ u, v_1, \ldots, v_{n-1} \}$ is a positively
oriented basis for $T_p\m$, whenever $u\in T_p\m$ is outward pointing. We
will denote by $\nu$ the outward pointing unit conormal vector field along
$\partial\m$. By means of the diffeomorphism
$\psi|_{\partial M}:\partial\m\fle\Sigma$, the orientation of $\partial\m$ is
induced on each connected component of $\Sigma$. On each connected component $P_0$ of $P$, we distinguish a connected component $\Sigma_0\subset P_0$ of $\Sigma$. 
Let $\eta_0$ be the
unitary vector field normal to $\Sigma_0$ in $P_0$ which points outward with respect to the domain
in $P_0$ bounded by $\Sigma_0$. Now, we choose $\xi_0$ the unique unitary vector field normal to $P_0$ in \mbar\ 
which is compatible with $\eta_0$ and with the orientation of $\Sigma_0$. We note that the chosen orientation of $P_0$ given by the field $\xi_0$ determines a unique choice to the
unitary vector field $\eta$ normal to each components of $\Sigma$ in $P_0$ such that $\eta|_{\Sigma_0}=\eta_0$. We repeat this process to the others connected components of $P$ and hence we obtain unitary vector fields $\eta$ normal to $\Sigma$ in $P$, and $\xi$ normal to $P$ in \mbar\ . With
this choice, given a point $p\in\Sigma$, a basis $\{ v_1, \ldots, v_{n-1} \}$ for $T_p\Sigma$ is
positively oriented if and only if $\{ \eta(p), v_1, \ldots, v_{n-1} \}$ is a positively
oriented basis for $T_pP$.

Let $\{ e_1, \ldots, e_{n-1} \}$ be a (locally defined) positively
oriented frame field along a fixed connected component of $\partial\m$. Using this frame, we can write
$\nu=e_1\times\ldots\times e_{n-1}\times {\bf N}$, and similarly
$\eta=e_1\times\ldots\times e_{n-1}\times\xi$, since
$\det(\nu,e_1,\ldots,e_{n-1},{\bf N})=1=\det(\eta,e_1,\ldots,e_{n-1},\xi)$.
From these expressions we easily compute
\begin{eqnarray*}
\eta & = & e_1\times\ldots\times e_{n-1}\times\xi= e_1\times\ldots\times 
e_{n-1}\times(\g{\xi}{{\bf N}}{\bf N}+\g{\xi}{\nu}\nu)\\
{} & = &   \g{\xi}{{\bf N}}\nu-\g{\xi}{\nu}{\bf N}, 
\end{eqnarray*}
that is,
\beq
\label{nueta}
\g{\eta}{\nu}=\g{\xi}{{\bf N}} \quad \mathrm{and} \quad
\g{\eta}{{\bf N}}=-\g{\xi}{\nu}.
\eeq

Let $A_\Sigma$ (respectively, $A_P$) denote the shape operator of
$\Sigma^{n-1}\subset P^n$ (respectively, $P^n\subset\mbar^{n+1}$) with
respect to the unit normal vector field $\eta$ (respectively, $\xi$). It
then follows that
\[
\nablabar_{e_i}e_j=\sum_{k=1}^{n-1}\g{\nablabar_{e_i}e_j}{e_k}e_k+
\g{\nablabar_{e_i}e_j}{\nu}\nu+\g{Ae_i}{e_j}{\bf N},
\]
for every $1\leq i,j\leq n-1$, and also
\[
\nablabar_{e_i}e_j=\sum_{k=1}^{n-1}\g{\nablabar_{e_i}e_j}{e_k}e_k+
\g{A_{\Sigma}e_i}{e_j}\eta+\g{A_Pe_i}{e_j}\xi,
\]
so that from \rf{nueta} we have that
\beq
\label{formaA1}
\g{Ae_i}{e_j}=-\g{A_{\Sigma}e_i}{e_j}\g{\xi}{\nu}+\g{A_Pe_i}{e_j}\g{\xi}{{\bf N}}.
\eeq

Equality \rf{formaA1} above shows us that it is not possible to go further
without any additional geometric hypothesis on the geometry of the
inclusion $P\subset\mbar$. A hypothesis of relevant geometric nature, and
which is also technically quite appropriate for us, consists on
assuming the umbilicity of $P\subset\mbar$. Then, from now on let us
suppose that $P$ is a totally umbilical hypersurface in \mbar. Therefore,
there exists a smooth function $\lambda\in\mathcal{C}^\infty(P)$ such that
$A_P=\lambda I$, where $I$ denotes the identity in $\mathcal{X}(P)$, and \rf{formaA1} becomes
\beq
\label{formaA2}
\g{Ae_i}{e_j}=-\g{A_{\Sigma}e_i}{e_j}\g{\xi}{\nu}+\lambda\g{\xi}{{\bf N}}\delta_{ij},
\quad 1\leq i,j\leq n-1.
\eeq

We now suppose that the basis $\{ e_1, \ldots, e_{n-1} \}\subset T_p(\partial M)$ on the boundary
is chosen such that it is formed by eigenvectors of $A_\Sigma$, and
let us denote its corresponding eigenvalues by $\tau_1(p), \ldots, \tau_{n-1}(p)$. In other words,
\[
A_\Sigma e_i=\tau_ie_i, \quad 1\leq i\leq n-1.
\]
Hence by \rf{formaA2}, $\g{Ae_i}{e_j}=0$ when $i\neq j$, and for each
$p\in\partial\m$, the matrix of $A$ in the orthonormal basis $\{ e_1, \ldots, e_{n-1}, \nu \}$
of $T_pM$ is given by
\beq
\label{matrixA}
A=\left(
\begin{array}{ccccc}
\gamma _1  & 0 & \cdots  & 0 & \left\langle A\nu ,e_1 \right\rangle  \\
0 & \gamma _2  & \cdots  & 0 & \left\langle A\nu ,e_2 \right\rangle  \\
\vdots  & \vdots  & \ddots & \vdots & \vdots  \\
0 & 0 & \cdots  & \gamma _{n-1} & \left\langle A\nu ,e_{n-1} \right\rangle \\
\left\langle A\nu ,e_1 \right\rangle  & \left\langle A\nu ,e_2 \right\rangle  & \cdots  &
\left\langle A\nu ,e_{n-1}\right\rangle  & \left\langle A\nu,\nu\right\rangle
\end{array}
\right),
\eeq
where $\gamma_i=-\tau_i\g{\xi}{\nu}+\lambda\g{\xi}{{\bf N}}$ for $1\leq i\leq
n-1$.

Now we compute the characteristic polynomial of $A$. To do that, we begin
by observing that
\begin{eqnarray}
\label{det}
\det (tI_n-A) & = & (t-\gamma_{n-1})\det(tI_{n-1}-\Lambda(\gamma_1, \ldots, \gamma_{n-2}))\\
\nonumber {} & - & \g{A\nu}{e_{n-1}}^2(t-\gamma_1)\ldots (t-\gamma_{n-2})
\end{eqnarray}
where
\[
\Lambda(\gamma_1, \ldots, \gamma_{n-2})=
\left(
\begin{array}{ccccc}
\gamma _1  & 0 & \cdots  & 0 & \left\langle A\nu ,e_1 \right\rangle  \\
0 & \gamma _2 & \cdots  & 0 & \left\langle A\nu ,e_2\right\rangle \\
\vdots  & \vdots  & \ddots & \vdots & \vdots  \\
0 & 0 & \cdots  & \gamma _{n-2}  & \left\langle A\nu ,e_{n-2} \right\rangle \\
\left\langle A\nu ,e_1 \right\rangle  & \left\langle A\nu ,e_2 \right\rangle  & \cdots  &
\left\langle A\nu ,e_{n-2}\right\rangle & \left\langle A\nu,\nu\right\rangle
\end{array}
\right).
\]
Therefore, applying a simple induction argument on $n$ in \rf{det}, we
obtain that the characteristic polynomial of $A$ is given by
\begin{eqnarray*}
\det(tI_n-A) & = & (t-\langle A\nu,\nu\rangle)\sum_{i=0}^{n-1}(-1)^{i}s_i(\gamma)t^{n-1-i} \\
{} & - & \sum_{i=1}^{n-1}\langle A\nu,e_i\rangle^2
\sum_{j=0}^{n-2}(-1)^{j}s_j(\widehat{\gamma_i})t^{n-2-j},
\end{eqnarray*}
where $s_r(\gamma)$ (respectively $s_r(\widehat{\gamma_i})$) stands for the elementary symmetric
functions of $\gamma_1,\ldots,\gamma_{n-1}$, (respectively
$\gamma_1,\ldots,\widehat{\gamma_i},\ldots,\gamma_{n-1}$), and, as usual,
$s_0(\gamma)=s_0(\widehat{\gamma_i})=1$ by definition.
Comparing the terms of above polynomials, we conclude from \rf{poly} that
the symmetric function of curvature $S_r$ of the hypersurface \m, at a boundary point
$p\in\partial\m$, is given by
\begin{eqnarray}
S_1&=&s_1(\gamma)+\g{A\nu}{\nu},\label{S1}\\
S_2&=&s_2(\gamma)+s_1(\gamma)\g{A\nu}{\nu}-\sum_{i=1}^{n-1}\g{A\nu}{e_i}^2,\label{S2}\\
S_r&=&s_r(\gamma)+s_{r-1}(\gamma)\g{A\nu}{\nu}-
\sum_{i=1}^{n-1}s_{r-2}(\widehat{\gamma_i})\g{A\nu}{e_i}^2, \label{Sr}
\end{eqnarray}
for $3\leq r\leq n$. 

\section{The Newton transformations on the boundary}
Observe that expressions \rf{S1} \rf{S2} and \rf{Sr} provide us with a partial 
answer to our initial question, since it relates the geometry of the 
hypersurface \m\ along its boundary $\partial\m$ (given by the 
$r$-curvature $S_r$) to the geometry of $\Sigma\subset P$ and the geometry 
of $P\subset\mbar$ (given by $s_r(\gamma)$). But this expression it is not 
still satisfactory for our purposes. We need the following essential 
auxiliary result. 
\begin{lemma}
\label{lemaPrnu}
Let $P^n\subset\mbar$ be an orientable totally umbilical hypersurface in \mbar\, and let
$\Sigma\subset P$
be an orientable $(n-1)$-dimensional compact submanifold in $P^n$. Let \x\ be an orientable 
connected
hypersurface with boundary $\Sigma=\psi(\partial\m)$, and let $\nu$ stands for the outward
pointing unit conormal vector field along $\partial\m\subset\m$. Then, along the boundary
$\partial\m$ and for every $1\leq r\leq n-1$, it holds
\beq
\label{Prnu}
\g{T_r\nu}{\nu}=s_r(\gamma)=s_r(\gamma_1, \ldots, \gamma_{n-1}),
\eeq
where $\gamma_i=-\tau_i\g{\xi}{\nu}+\lambda\g{\xi}{{\bf N}}$ for $1\leq i\leq
n-1$. Here $\tau_1, \ldots, \tau_{n-1}$ are the principal curvatures of
$\Sigma\subset P$ with respect to the outward pointing unitary
normal, ${\bf N}$ is the unitary normal field of \m, $\xi$ is the unitary normal field of
$P\subset\mbar$, and $\lambda$ is the umbilicity factor of $P\subset\mbar$ (with respect to
$\xi$).
\end{lemma}
\begin{proof}
We will use induction on $r$. First, observe that from (\ref{S1}) it follows that \rf{Prnu} holds for
$r=1$. For a given $2\leq r\leq n-1$, suppose that
\beq
\label{induction}
\g{T_j\nu}{\nu}=s_j(\gamma)
\eeq
holds for all $1\leq j\leq r-1$. Observe that
\[
A\nu=\sum_{i=1}^{n-1}\g{A\nu}{e_i}e_i+\g{A\nu}{\nu}\nu,
\]
so that from the inductive definition of $T_r$ and \rf{induction} we conclude that
\begin{eqnarray}
\label{Prnu2}
\g{T_r\nu}{\nu} & = & S_r-\g{T_{r-1}\nu}{A\nu} \\
\nonumber {} & = &
S_r-\g{T_{r-1}\nu}{\nu}\g{A\nu}{\nu}-\sum_{i=1}^{n-1}\g{T_{r-1}\nu}{e_i}\g{A\nu}{e_i}\\
\nonumber {} & = &
S_r-s_{r-1}(\gamma)\g{A\nu}{\nu}-\sum_{i=1}^{n-1}\g{T_{r-1}\nu}{e_i}\g{A\nu}{e_i}.
\end{eqnarray}
On the other hand, we also know from \rf{matrixA} that
\[
Ae_i=\gamma_ie_i+\g{A\nu}{e_i}\nu,
\]
so that from our induction hypothesis \rf{induction} we have for every
$1\leq j\leq r-1$,
\[
\g{T_j\nu}{e_i}=-\g{T_{j-1}\nu}{Ae_i}=
-\gamma_i\g{T_{j-1}\nu}{e_i}-s_{j-1}(\gamma)\g{A\nu}{e_i}.
\]
This implies by a recursive argument that
\beq
\label{Pr-1}
\g{T_{r-1}\nu}{e_i}=
-\g{A\nu}{e_i}\sum_{j=0}^{r-2}(-1)^js_{r-2-j}(\gamma)\gamma_i^j
=-\g{A\nu}{e_i}s_{r-2}(\widehat{\gamma_i}),
\eeq
since it is not difficult to see that
\[
s_m(\widehat{\gamma_i})=\sum_{j=0}^{m}(-1)^js_{m-j}(\gamma)\gamma_i^j
\]
for every $1\leq m\leq n-1$. Using now \rf{Pr-1} in \rf{Prnu2}, along with
\rf{Sr}, we conclude that
\[
\g{T_r\nu}{\nu}=S_r-s_{r-1}(\gamma)\g{A\nu}{\nu}+\sum_{i=1}^{n-1}s_{r-2}(\widehat{\gamma_i})
\g{A\nu}{e_i}^2=s_r(\gamma).
\]
This finishes the proof of Lemma\rl{lemaPrnu}.
\end{proof}

Now, it remains to know how the elementary symmetric function $s_r(\gamma)$ can be
expressed in terms of the principal curvatures $\tau_1, \ldots,
\tau_{n-1}$ of the inclusion $\Sigma\subset P$ and the umbilicity factor $\lambda$
of $P\subset\mbar$. To see this, let us write $\gamma_i=\alpha_i+\beta$,
where $\alpha_i=-\tau_i\g{\xi}{\nu}$ and $\beta=\lambda\g{\xi}{{\bf N}}$, for
each $i=1, \ldots, n-1$.
\begin{lemma}
\[
s_r(\gamma)=\sum_{j=0}^{r}{{n-1-j}\choose{r-j}}\beta^{r-j}s_j(\alpha),\qquad 1\leq
r\leq n-1.
\]
\end{lemma}
\begin{proof}
Recall that $s_r(\gamma)$ can be defined by the following polynomial
identity \rf{poly},
\[
\sum_{r=0}^{n-1}(-1)^rs_r(\gamma)t^{n-1-r}=(t-\gamma_1)\cdots(t-\gamma_{n-1}).
\]
Since each $\gamma_i=\alpha_i+\beta$, the right hand side of this equality
can be written as follows
\[
((t-\beta)-\alpha_1)\cdots((t-\beta)-\alpha_{n-1})=
\sum_{j=0}^{n-1}(-1)^{j}s_{j}(\alpha)(t-\beta)^{n-1-j}.
\]
On the other hand, computing the right hand side of this last equality, we
obtain
\[
\sum_{j=0}^{n-1}(-1)^{j}s_{j}(\alpha)(t-\beta)^{n-1-j}=
\sum_{j=0}^{n-1}\sum_{k=0}^{n-1-j}(-1)^{k+j}{{n-1-j}\choose{k}}
\beta^k s_{j}(\alpha)t^{n-1-k-j},
\]
which after a re-ordering becomes
\[
\sum_{r=0}^{n-1}(-1)^r
\left(\sum_{j=0}^{r}{{n-1-j}\choose{r-j}}\beta^{r-j}s_{j}(\alpha)\right)
t^{n-1-r}.
\]
Therefore, we have obtained the following equality between polynomials
\[
\sum_{r=0}^{n-1}(-1)^rs_r(\gamma)t^{n-1-r}=
\sum_{r=0}^{n-1}(-1)^r
\left(\sum_{j=0}^{r}{{n-1-j}\choose{r-j}}\beta^{r-j}s_{j}(\alpha)\right)
t^{n-1-r},
\]
which concludes the proof.
\end{proof}

We summarize what was made above as follows.
\begin{proposition}
\label{propositionPrnu}
Let $P^n\subset\mbar$ be an orientable totally umbilical hypersurface in \mbar\, and let
$\Sigma\subset P$
be an orientable $(n-1)$-dimensional compact submanifold in $P^n$. Let \x\ be an orientable
hypersurface with boundary $\Sigma=\psi(\partial\m)$, and let $\nu$ stands for the outward
pointing unit conormal vector field along $\partial\m\subset\m$. Then, along the boundary
$\partial\m$ and for every $1\leq r\leq n-1$, it holds that
\beq
\label{Prnu3}
\g{T_r\nu}{\nu}=\sum_{j=0}^{r}(-1)^j{{n-1-j}\choose{r-j}}\lambda^{r-j}\g{\xi}{{\bf N}}^{r-j}
\g{\xi}{\nu}^{j}s_j.
\eeq
Here $s_j=s_j(\tau_1,\ldots,\tau_{n-1})$, $0\leq j\leq n-1$, are the elementary symmetric
functions of $\tau_1, \ldots, \tau_{n-1}$, the principal
curvatures of $\Sigma\subset P$ with respect to the outward pointing unitary
normal, ${\bf N}$ is the unitary normal field of \m, $\xi$ is the unitary normal field of
$P\subset\mbar$, and $\lambda$ is the umbilicity factor of $P\subset\mbar$ (with respect to
$\xi$).
\end{proposition}

\section{Transversality versus Ellipticity}
The relationship between the $S_r$'s and the $s_r(\gamma)$'s given in 
\rf{S1}, \rf{S2} and \rf{Sr}, as well as the expression for $\g{T_r\nu}{\nu}$ given in 
\rf{Prnu3} becomes specially simple in the case where the inclusion $P\subset\mbar$ is 
\textit{totally geodesic}, that is, when $\lambda=0$. In that case 
$\gamma_i=-\tau_i\g{\xi}{\nu}$, and we have the following.
\begin{corollary}
\label{c1}
Let $\Sigma$ be an orientable $(n-1)$-dimensional compact submanifold in
an orientable totally geodesic hypersurface $P^n\subset\mbar^{n+1}$. Let \x\ be an
orientable hypersurface with boundary $\Sigma=\psi(\partial\m)$, and let $\nu$ stands for the
outward pointing unit conormal vector field along $\partial\m\subset\m$. Then, along the
boundary $\partial\m$ and for every $1\leq r\leq n$, it holds that
\begin{eqnarray}
S_1 & = & -s_1\g{\xi}{\nu}+\g{A\nu}{\nu},\label{S1bis}\\
S_2 & = & s_2\g{\xi}{\nu}^{2}-s_1\g{\xi}{\nu}\g{A\nu}{\nu}-
\sum_{i=1}^{n-1}\g{A\nu}{e_i}^2,\label{S2bis}\\
S_r & = & (-1)^rs_r\g{\xi}{\nu}^{r}+(-1)^{r-1}s_{r-1}\g{\xi}{\nu}^{r-1}\g{A\nu}{\nu}\label{Srbis}\\
\nonumber {} & {} & 
-(-1)^{r-2}\g{\xi}{\nu}^{r-2}\sum_{i=1}^{n-1}s_{r-2}(\widehat{\tau_i})\g{A\nu}{e_i}^2, 
\end{eqnarray}
for $3\leq r\leq n$, and
\beq
\label{Prnu4}
\g{T_r\nu}{\nu}=(-1)^rs_r\g{\xi}{\nu}^{r},
\eeq
where $s_n=0$ and for every $1\leq r\leq n-1$
\[
s_r=s_r(\tau_1,\ldots,\tau_{n-1})
\]
is the $r$-th elementary symmetric function of $\tau_1,\ldots,\tau_{n-1}$, the principal
curvatures of $\Sigma\subset P$ with respect to the outward pointing unitary
normal, and $\xi$ is the unitary normal field of $P\subset\mbar$.
\end{corollary}

It is not difficult to see that \rf{Prnu4} establishes a very strong
relationship between the \textit{transversality} of \m\ with respect to
$P$ along the boundary $\partial\m$, and the \textit{ellipticity} on \m\ of
the $r$-th Newton transformation $T_r$, when $r\geq 1$ (recall that $T_0=I$).
That relationship between transversality and ellipticity will be actually one of the 
keys of the proof of our symmetry results (Theorem\rl{calotaRn}, Theorem\rl{simetriaHn} 
and Theorem\rl{simetriaSn}).
In fact, saying that \m\ is not transverse to $P$ along its boundary $\partial\m$ means that 
there exists a point $p\in\partial\m$ such that $\g{\xi}{\nu}(p)=0$, which implies from 
\rf{Prnu4} that $\g{T_r\nu}{\nu}(p)=0$, $r\geq 1$. Therefore we can conclude that if the Newton 
transformation $T_r$ is positive definite on \m\ for some $1\leq r\leq n-1$, then the 
hypersurface \m\ is necessarily transverse to $P$ along its boundary. 

Observe that in the case where $S_n$ does not vanish on \m\ and $n\geq 3$, transversality easily follows from
expression \rf{Srbis}. In fact, by \rf{Srbis} we have along the boundary $\partial\m$
\begin{eqnarray*}
S_n & = & 
(-1)^{n-1}s_{n-1}\g{\xi}{\nu}^{n-1}\g{A\nu}{\nu}\\
{} & {} & 
+(-1)^{n-1}\g{\xi}{\nu}^{n-2}\sum_{i=1}^{n-1}s_{n-2}(\widehat{\tau_i})\g{A\nu}{e_i}^2.
\end{eqnarray*}
In particular, if there exists a point $p\in\partial\m$ where $\g{\xi}{\nu}(p)=0$, then
$S_n(p)=0$ (since $n\geq 3$). In the same way, if we assume that $n\geq 2$ 
and $S_2$ is positive everywhere on \m, then \rf{S2bis} also implies that \m\ is transverse to 
$P$ along the boundary.

We summarize what was made above as follows.
\begin{proposition}
\label{transverse}
Let $\Sigma$ be an orientable $(n-1)$-dimensional compact submanifold in
an orientable totally geodesic hypersurface $P^n\subset\mbar^{n+1}$ and let \x\ be an
orientable hypersurface with boundary $\Sigma=\psi(\partial\m)$. Then each 
one of the following hypothesis individually implies that \m\ is transverse to $P$ 
along the boundary $\partial\m$:
\begin{itemize}
\item For a given $1\leq r\leq n-1$, the Newton transformation $T_r$ is definite positive on 
\m.
\item $n\geq 3$ and $S_n\neq 0$ on \m.
\item $S_2>0$ on \m.
\end{itemize}
\end{proposition}

\section{Symmetry for hypersurfaces in Euclidean space}
The totally umbilic hypersurfaces of Euclidean space \Rn\ are the totally 
geodesic hyperplanes and the round $n$-spheres. They trivially have constant 
$r$-mean curvature for each $r=0,\ldots, n$. Actually, the hyperplanes have vanishing 
$r$-mean curvature $H_r=0$, and, after an appropriate choice of the unit 
normal vector field, the round $n$-spheres of radius $\varrho>0$ have 
constant $r$-mean curvature $H_r=1/\varrho^r$. Let us fix a hyperplane 
$\Pi\subset\Rn$ and an $(n-1)$-sphere $\Sigma\subset\Pi$. Then the hyperplanar round ball 
bounded by $\Sigma$ in $\Pi$, and the spherical caps bounded by $\Sigma$ (of radii greater than 
or equal to the radius of $\Sigma$) are examples of compact hypersurfaces 
embedded into \Rn\ with constant $r$-mean curvature and bounded by $\Sigma$. 
In this context, it was conjectured in \cite{BEMR} that these examples are 
the only compact embedded hypersurfaces in \Rn\ with constant mean 
curvature and spherical boundary. 
Related to this conjecture we have the 
following symmetry theorem for hypersurfaces in Euclidean space \cite{AM}.
\begin{theorem}
\label{calotaRn}
Let $\Sigma$ be an strictly convex compact $(n-1)$-dimensional submanifold in a
hyperplane $\Pi\subset\Rn$, and let \xR\ be a compact embedded hypersurface with boundary
$\Sigma$. Let us assume that for a given $2\leq r\leq n$, the $r$-mean
curvature $H_r$ of \m\ is a nonzero constant . Then \m\ has all the
symmetries of $\Sigma$. In particular, if the boundary $\Sigma$ is a round
$(n-1)$-sphere of \Rn, then \m\ is a spherical cap.
\end{theorem}
\begin{proof}
It is not difficult to see that under the hypothesis above there exists at
least one interior elliptic point of \m, that is, an interior point of
\m\ where, after an appropriate orientation of \m, all the principal curvatures are positive. 
In fact, since \m\ is not part of a hyperplane (because of $H_r\neq 0$), then one easily finds 
a radius $R>0$ and a point $a\in\Rn$ such that the closed round ball $\overline{B}(a,R)$
contains \m\ and such that there is a point $p_0\in\mathrm{int}(\m)\cap\partial B(a,R)$
(englobe \m\ with spheres of large radius until such a sphere touches \m\ on one side at an
interior point). In particular, in the chosen orientation the constant $H_r=H_r(p_0)>0$ is 
positive. The existence of an elliptic point, jointly with the fact that $H_r$ is a positive
constant, allows us to conclude that the Newton transformation $T_{r-1}$
is positive definite on \m\ (see \cite[Proposition 3.2]{BC} and \cite[p. 232]{Ro}).
Therefore, from Proposition\rl{transverse} it follows that \m\ is
transverse to $\Pi$ along the boundary $\partial\m$. Our result then is a consequence of
Theorem 7.3 in \cite{Ro}.
\end{proof}
As a consequence of Theorem\rl{calotaRn} we can conclude that the
conjecture of the spherical cap \cite{BEMR} is true for the case of embedded hypersurfaces with
constant $r$-mean curvature in \Rn, when $r\geq 2$ \cite{AM}.
\begin{corollary}
The only  compact embedded hypersurfaces in \Rn\ with constant $r$-mean curvature $H_r$
(with $2\leq r\leq n$) and spherical boundary are the hyperplanar round balls (with $H_r=0$) 
and the spherical caps (with $H_r$ a nonzero constant).
\end{corollary}
Indeed, if  \m\ is not a hyperplanar round ball, then the constant $r$-mean curvature
must be necessarily nonzero because there exists at least one interior elliptic point of \m. 
In particular, when $r=2$ saying that $H_2$ is constant is equivalent to 
saying that the scalar curvature is constant (see equation \rf{scalar}), 
so that the result reads as follows.
\begin{corollary}
The only compact  embedded hypersurfaces in \Rn\ with constant scalar curvature and spherical
boundary are the hyperplanar round balls (with zero scalar curvature) and the spherical caps
(with positive constant scalar curvature).
\end{corollary}

Our objective in Sections\rl{hyperbolic} and \rl{sphere} is to extend the
symmetry result given in Theorem\rl{calotaRn} to the case of hypersurfaces in hyperbolic space 
and hypersurfaces in sphere, as well as the corresponding solution to the spherical cap 
conjecture for the case of constant $r$-mean curvature, $r\geq 2$. A result of this type was 
first given by Nelli and Rosenberg in \cite[Theorem 3.1]{NR} for hypersurfaces with constant
mean curvature in hyperbolic space. On the other hand, the corresponding 
result for the case of hypersurfaces with constant mean curvature in 
sphere has been recently given by Lira \cite{Lira}.
As observed by Nelli and Rosenberg,
their result could be extended to the case of constant $r$-mean curvature as
soon as a certain \textit{flux formula} could be established. In the next
section, we will derive such a flux formula.

\section{A flux formula}
\label{sectionflux}
In this section we will derive a general flux formula for the geometric
configuration considered in Section\rl{configuration} in the case where
the Riemannian ambient space \mbar\ is equipped with a conformal vector field
$Y\in\mathcal{X}(\mbar)$. Recall that the fact that $Y$ is conformal means that the Lie derivative
of the metric tensor of \mbar\ with respect to $Y$ satisfies
\[
\mathcal{L}_Y\langle , \rangle=2\phi\langle , \rangle
\]
for a certain smooth function $\phi\in\mathcal{C}^\infty(\mbar)$. In other words,
\begin{equation}
\label{conformal}
\g{\nablabar_VY}{W}+\g{V}{\nablabar_WY}=2\phi\g{V}{W},
\end{equation}
for every vector fields $V,W\in\mathcal{X}(\mbar)$.

In order to derive our general flux formula, let us consider $Y^\top\in\xm$ the vector
field obtained on the hypersurface \m\ by taking the tangential component
of $Y$, that is, $Y^\top=Y-f{\bf N}$, where $f=\langle Y,{\bf N}\rangle$.
Most of the interesting and useful integral formulas in Riemannian geometry are obtained by
computing the divergence of certain vector fields and applying the divergence theorem. The
interesting integral formulas therefore correspond to vector fields with interesting
divergences. Our idea here is to compute the divergence $\Div(T_rY^\top)$.
Using that $\nabla_UT_r$ is self-adjoint for any tangent vector field $U\in\xm$, an
easy computation shows that
\begin{equation}
\label{e1}
\Div(T_rY^\top)=\g{\Div T_r}{Y}+
\sum_{i=1}^n\g{\nabla_{e_i}Y^\top}{T_re_i},
\end{equation}
where $\{ e_1, \ldots, e_n \}$ is a local orthonormal frame on \m\ and
$\Div T_r$ is given by \rf{div} in Lemma\rl{divlema}. From the conformal equation
(\ref{conformal}), we obtain
\begin{eqnarray*}
2\phi\langle T_rU,U\rangle & = &
\g{\nablabar_{T_rU}Y}{U}+\g{\nablabar_U Y}{T_rU} \\
{} & = & \langle\nablabar_{T_rU}Y^\top, U\rangle+f\langle \nablabar_{T_rU}{\bf N},U\rangle+
\langle\nablabar_U Y^\top, T_rU\rangle +f\langle\nablabar_U {\bf N}, T_rU\rangle \\
{} & = & \langle\nabla_{T_rU}Y^\top, U\rangle+\langle \nabla_UY^\top,
T_rU\rangle-f\langle AT_rU, U\rangle-f\langle AU, T_rU\rangle,
\end{eqnarray*}
that is
\begin{equation}
\label{eq1}
\langle \nabla_{T_rU}Y^\top, U\rangle+\langle \nabla_UY^\top,
T_rU\rangle = 2\phi\langle T_rU,U\rangle+2f\langle AT_rU, U\rangle
\end{equation}
Let us choose $\{e_1, \ldots,e_n \}$ a local
orthonormal frame on \m\ such that diagonalizes $A$. We know then that it
also diagonalizes $T_r$ with eigenvalues $\mu_{1,r}, \ldots, \mu_{n,r}$, and
therefore
\[
\langle\nabla_{e_i}Y^\top,T_re_i\rangle=\mu_{i,r}\langle\nabla_{e_i}Y^\top,e_i\rangle=
\langle e_i,\nabla_{T_re_i}Y^\top\rangle,
\]
so that from (\ref{eq1}) we obtain
\[
\langle\nabla_{e_i}Y^\top,{T_re_i}\rangle=\phi\langle{e_i},{T_re_i}\rangle+
\langle{Y},{{\bf N}}\rangle\langle{AT_re_i},{e_i}\rangle.
\]
Taking trace here and using \rf{traces} and \rf{tracesbis}, equation (\ref{e1}) becomes
\begin{equation}
\label{divPr}
\Div(T_rY^\top)=\g{\Div T_r}{Y}+c_r(\phi H_r+\g{Y}{{\bf N}}H_{r+1}),
\end{equation}
where $c_r=(r+1)\binom{n}{r+1}$.
Integrating now \rf{divPr} on \m, the Stokes theorem implies the following
integral formula for every $0\leq r\leq n-1$,
\begin{eqnarray}
\label{phi1}
\nonumber \oint_{\partial M}\g{T_r\nu}{Y}ds & = & \int_M\Div(T_rY^\top)dM \\
{} & = & \int_M\g{\Div T_r}{Y}dM+c_r\int_M(\phi H_r+\g{Y}{{\bf N}}H_{r+1})dM
\end{eqnarray}
Here $dM$ denotes the $n$-dimensional volume element of \m\ with respect
to the induced metric and the chosen orientation, and $ds$ is the
$(n-1)$-dimensional volume element induced on $\partial\m$.

On the other hand, let $D^n$ be a compact orientable hypersurface in
\mbar\ with smooth boundary that satisfies $\partial D=\partial\m$, such that $M\cup D$ is an 
oriented $n$-cycle of \mbar, with $D$ oriented by the unit normal field $n_D$. We suppose that 
$\m\cup D=\partial\Omega$, where $\Omega$ is a
compact oriented domain immersed in \mbar. From the conformal equation
\rf{conformal}, we easily see that $\Divbar Y=(n+1)\phi$. Therefore, from
Stokes theorem we obtain that
\beq
\label{phi2}
\int_M\g{Y}{{\bf N}}dM=-\int_D\g{Y}{n_D}dD+(n+1)\int_\Omega\phi d\mbar,
\eeq
where $dD$ stands for the $n$-dimensional volume element of $D$ with
respect to the orientation given by $n_D$, and $d\mbar$ denotes the
$(n+1)$-dimensional volume element on \mbar.
Now, from \rf{phi1} and \rf{phi2} we conclude the following general
\textit{flux formula}.
\begin{proposition}
\label{propflux}
Let \x\ be an immersed compact orientable hypersurface with boundary $\partial\m$, and
let $D^n$ be a compact orientable hypersurface with boundary $\partial
D=\partial\m$. Assume that $\m\cup D$ is an oriented $n$-cycle of \mbar,
and let ${\bf N}$ and $n_D$ be the unit normal fields which orient \m\ and $D$, respectively.
If the $r$-mean curvature $H_r$ is constant, $1\leq r\leq n$, then for every conformal vector
field $Y\in\mathcal{X}(\mbar)$ the following formula holds
\begin{eqnarray}
\label{fluxo}
\nonumber \oint_{\partial M}\g{T_{r-1}\nu}{Y}ds & = & \int_M\g{\Div T_{r-1}}{Y}dM+
r\binom{n}{r}\int_M\phi H_{r-1}dM \\
{} & {} & -r\binom{n}{r}H_r\int_D\g{Y}{n_D}dD+(n+1)r\binom{n}{r}H_r\int_\Omega\phi d\mbar,
\end{eqnarray}
$\nu$ is the outward pointing conormal to \m\ along $\partial\m$.
\end{proposition}

Formula \rf{fluxo} becomes specially simple when the ambient space \mbar\
has constant sectional curvature, and the field $Y$ is a Killing vector field,
that is, $\phi=0$. In that case, the Newton transformations are
divergence-free (Corollary\rl{divcoro}) and from formula \rf{fluxo} we
derive the \textit{balancing formula} given by Rosenberg in \cite[Theorem 7.2]{Ro}
(see also \cite{Ku,BE1,BEMR,NR} for the case of constant mean curvature).
\begin{corollary}
\label{coroflux}
If \mbar\ has constant sectional curvature, then for every Killing vector
field $Y\in\mathcal{X}(\mbar)$ the flux formula becomes
\begin{equation}
\label{fluxo3}
\oint_{\partial M}\g{T_{r-1}\nu}{Y}ds=-r\binom{n}{r}H_r\int_D\g{Y}{n_D}dD,
\end{equation}
where $\nu$ is the outward pointing conormal to \m\ along $\partial\m$.
\end{corollary}

On the other hand, when the ambient space \mbar\ has constant sectional curvature, and the field
$Y$ is a homothetic (and non-Killing) vector field, then we may assume without loss of generality
that $\phi=1$ and \rf{fluxo} becomes
\begin{eqnarray}
\label{fluxo2}
\nonumber \oint_{\partial M}\g{T_{r-1}\nu}{Y}ds & = & -r\binom{n}{r}H_r\int_D\g{Y}{n_D}dD \\
{} & {} & +r\binom{n}{r}\int_M H_{r-1}dM+(n+1)r\binom{n}{r}H_r\mathrm{vol}(\Omega).
\end{eqnarray}
As a consequence of \rf{fluxo2} we obtain the following flux formula for
$r$-minimal hypersurfaces.
\begin{proposition}
Let \x\ be a compact orientable hypersurface with boundary $\partial\m$
immersed into a Riemannian space of constant sectional curvature. If \m\
is $r$-minimal in \mbar, that is, $H_r=0$, then for every homothetic (non-Killing)
vector field $Y\in\mathcal{X}(\mbar)$ the following formula holds
\begin{eqnarray}
\label{fluxo4}
\oint_{\partial M}\g{T_{r-1}\nu}{Y}ds=r\binom{n}{r}\int_M H_{r-1}dM.
\end{eqnarray}
\end{proposition}
In particular, for minimal hypersurfaces in Euclidean space with boundary 
in a round sphere we have the following consequence.
\begin{corollary}
\label{coroflux1}
Let $\Sigma$ be an orientable $(n-1)$-dimensional compact submanifold in a round 
sphere $\mathbb{S}^n(\varrho)\subset\Rn$ of radius $\varrho$, and let \xR\ be an immersed 
orientable compact minimal hypersurface with boundary 
$\Sigma=\psi(\partial\m)\subset\mathbb{S}^n(\varrho)$. Then
\[
\mathrm{vol}(\m)\leq\frac{\varrho}{n}\mathrm{vol}(\partial\m),
\]
and equality holds if and only if \m\ is orthogonal to $\mathbb{S}^n(\varrho)$
along the boundary $\partial\m$.
\end{corollary}
\begin{proof}
Consider the radial vector field $Y(p)=p$ in $\mathbb{R}^{n+1}$, which is a homothetic vector
field in \Rn\ with $\phi=1$, and let $\xi$ be the unit vector normal to $\mathbb{S}^n(\varrho)$.
Then, along $\mathbb{S}^n(\varrho)$ we have $Y=\varrho \xi$ and \rf{fluxo4} gives
\[
n\int_MdM=n\mathrm{vol}(M)=\oint_{\partial M}\g{\nu}{\varrho \xi}ds\leq\varrho
\oint_{\partial M}ds=\varrho \mathrm{vol}(\partial\m).
\]
Besides, equality holds if and only if $\xi=\nu$ along the boundary
$\partial\m$, or equivalently (see \rf{nueta}) $\g{{\bf N}}{\xi}=0$ along
$\partial\m$.
\end{proof}

Let us consider now the case of a hypersurface immersed into the hyperbolic space
\Hn. In that case, it will be appropriate to use the Minkowski space model
of hyperbolic space. Write $\mathbb{R}^{n+2}_1$ for $\mathbb{R}^{n+2}$
with the Lorentzian metric
\[
\g{}{}_1=-dx_0^2+dx_1^2+\cdots+dx_{n+1}^2.
\]
Then
\[
\Hn=\{ x\in\mathbb{R}^{n+2}_1 : \g{x}{x}_1=-1, \quad x_0>0 \}
\]
is a complete spacelike hypersurface in $\mathbb{R}^{n+2}_1$ with constant
sectional curvature $-1$ which provides the Minkowski space model for the
hyperbolic space.

Let $\Sigma\subset\Hn$ be an orientable $(n-1)$-dimensional compact submanifold in a
geodesic sphere $S(a,\varrho)$ of \Hn\ of center $a\in\Hn$ and geodesic radius 
$\varrho$, and let \xH\ be an orientable compact hypersurface with boundary 
$\Sigma=\psi(\partial M)$. 

Consider the vector field in \Hn\ represented in this model as 
$Y(p)=-a-\g{a}{p}p$ for every $p\in\Hn$. Observe that $Y$ is a conformal vector field in \Hn\ 
which is orthogonal to the geodesic spheres centered at the point $a$, with 
$\phi(p)=-\g{a}{p}=\cosh{(\tilde{\varrho}(p))}$ and 
$|Y(p)|=\sinh{(\tilde{\varrho}(p))}$, where 
$\tilde{\varrho}(p)=\mathrm{distance}(p,a)$ for every $p\in\Hn$.  
Therefore, along $S(a,\varrho)$ we have $Y=\sinh{\varrho}\xi$. 
Assume now that \m\ is minimal in \Hn. Then, it follows from \rf{fluxo} that 
\begin{eqnarray*}
\oint_{\partial M}\g{\nu}{Y}ds=\sinh{\varrho}\oint_{\partial M}\g{\nu}{\xi}ds=
n\int_M\cosh{(\tilde{\varrho})} dM.
\end{eqnarray*}
Thus, since $\cosh{(\tilde{\varrho})}\geq 1$, we conclude from here that
\begin{eqnarray*}
n\mathrm{vol}(M)\leq n\int_M\cosh{(\tilde{\varrho})} dM=
\sinh{\varrho}\oint_{\partial M}\g{\nu}{\xi}ds\leq
\sinh{\varrho}\mathrm{vol}(\partial M).
\end{eqnarray*}
Summing up, we have obtained the following result.
\begin{corollary}
\label{coroflux2} 
Let $\Sigma$ be an orientable $(n-1)$-dimensional compact submanifold in a
geodesic sphere $S(a,\varrho)$ of \Hn\ of center $a\in\Hn$ and geodesic radius 
$\varrho$, and let \xH\ be an immersed orientable compact minimal hypersurface with boundary 
$\Sigma=\psi(\partial M)\subset S(a,\varrho)$. Then
\[
\mathrm{vol}(M)\leq\frac{\sinh{\varrho}}{n}\mathrm{vol}(\partial M).
\]
\end{corollary}

Finally, let us consider the case of a hypersurface immersed into the sphere \Sn,
\[
\Sn=\{ x=(x_0,\ldots,x_{n+1})\in\mathbb{R}^{n+2} : \g{x}{x}=1 \}.
\]
Let $\Sigma$ be an orientable $(n-1)$-dimensional compact submanifold in a 
geodesic sphere $S(a,\varrho)$ of \Sn\ of center $a\in\Sn$ and radius $\varrho<\pi/2$, and let 
\xS\ be an orientable compact hypersurface with boundary 
$\Sigma=\psi(\partial\m)\subset S(a,\varrho)$.

In this case, consider the vector field in \Sn\ given by  
$Y(p)=-a+\g{a}{p}p$ for every $p\in\Sn$, with singularities at the focal points $\{a,-a\}$. 
Observe that $Y$ is a conformal vector field in \Sn\ 
which is orthogonal to the geodesic spheres centered at the point $a$, with 
$\phi(p)=\g{a}{p}=\cos{(\tilde{\varrho}(p))}$ and 
$|Y(p)|=\sin{(\tilde{\varrho}(p))}$, where 
$\tilde{\varrho}(p)=\mathrm{distance}(p,a)$ for every $p\in\Sn$.  
Therefore, along $S(a,\varrho)$ we have $Y=\sin{\varrho}\xi$. 
Assume now that \m\ is minimal in \Sn. Then, it follows from \rf{fluxo} that 
\begin{eqnarray*}
\oint_{\partial M}\g{\nu}{Y}ds=\sin{\varrho}\oint_{\partial M}\g{\nu}{\xi}ds=
n\int_M\cos{(\tilde{\varrho})} dM.
\end{eqnarray*}
Let us assume now that $M$ is contained in the open hemisphere centered at 
$a$. In that case, it is clear that 
$\min_M \cos{(\tilde{\varrho})}=\cos{\varrho_0}$,
where $\varrho_0=\max{\rm dist}(a,M)$, so that 
\begin{eqnarray*}
n\cos{\varrho_0}\mathrm{vol}(M)\leq n\int_M\cos{(\tilde{\varrho})} dM=
\sin{\varrho}\oint_{\partial M}\g{\nu}{\xi}ds\leq
\sin{\varrho}\mathrm{vol}(\partial M).
\end{eqnarray*}
This leads to the following result.
\begin{corollary} 
\label{coroflux3}
Let $\Sigma$ be an orientable $(n-1)$-dimensional compact submanifold in a
geodesic sphere $S(a,\varrho)$ of \Sn\ of center $a\in\Hn$ and geodesic radius 
$\varrho$, and let \xS\ be an immersed orientable compact minimal hypersurface with boundary 
$\Sigma=\psi(\partial M)\subset S(a,\varrho)$. Assume that $M$ is contained in the open 
hemisphere centered at $a$. Then 
\[
\mathrm{vol}(M)\leq\frac{\sin{\varrho}}{n\cos{\varrho_0}}\mathrm{vol}(\partial M),
\]
where $\varrho_0=\max{\rm dist}(a,M)$.
\end{corollary}

\section{Estimating the $r$-mean curvature by the geometry of the boundary}
\label{estimating}
In this section, we will describe an interesting application of our flux formula
\rf{fluxo3} and the formula \rf{Prnu4}. Let us consider the geometric configuration given
in Proposition\rl{transverse}; that is,
let $\Sigma$ be an orientable $(n-1)$-dimensional compact submanifold in an orientable
totally geodesic hypersurface $P^n\subset\mbar$, and let \x\ be an orientable compact
connected hypersurface with boundary $\Sigma=\psi(\partial\m)$ and constant $r$-mean curvature $H_r$.
Our objective here is to estimate $H_r$ by the geometry of the boundary. Assume that there exists a
Killing vector field $Y\in\mathcal{X}(\mbar)$ which is orthogonal to $P$. Then, we can write $Y$
along the boundary $\partial\m$ both as $Y=\g{Y}{\xi}\xi$ and also as
$Y=\g{Y}{\nu}\nu+\g{Y}{{\bf N}}{\bf N}$, and using \rf{Prnu4} we obtain
\[
\g{T_{r-1}\nu}{Y}=\g{Y}{\nu}\g{T_{r-1}\nu}{\nu}=(-1)^{r-1}s_{r-1}\g{Y}{\xi}\g{\xi}{\nu}^{r}
\]
along the boundary $\partial\m$.

Let us consider $D\subset P$ the domain in $P$ bounded by $\Sigma$, and let us orient $D$ by the
unit normal field $n_D$, so that $\m\cup D$ is an oriented $n$-cycle in \mbar.
Let us denote by $h_j$ the $j$-th mean curvature of $\Sigma\subset P$ with respect to the
outward pointing unitary normal $\eta$, that is,
\[
{{n-1}\choose{j}}h_j=s_j=s_j(\tau_1,\ldots,\tau_{n-1}), \quad 0\leq j\leq n-1.
\]
In the case where the ambient space \mbar\ has constant sectional
curvature, then our flux formula \rf{fluxo3} allows us to write
\beq
\label{estimativa}
nH_r\int_D\g{Y}{n_D}dD=(-1)^r\oint_{\partial M}h_{r-1}\g{Y}{\xi}\g{\xi}{\nu}^rds.
\eeq
Let us first apply formula \rf{estimativa} to the Euclidean case,  $\mbar=\Rn$.
\begin{theorem}
\label{estimate1}
Let $\Sigma$ be an orientable $(n-1)$-dimensional compact submanifold in a hyperplane
$P\subset\Rn$, and let \xR\ be an orientable immersed compact (connected)
hypersurface with boundary $\Sigma=\psi(\partial\m)$ and constant $r$-mean curvature $H_r$,
$1\leq r\leq n$. Then
\beq
\label{est1}
0\leq |H_r|\leq\frac{1}{n\ \mathrm{vol}(D)}\oint_{\partial M}|h_{r-1}|ds,
\eeq
where $h_{r-1}$ stands for the $(r-1)$-mean curvature of $\Sigma\subset
P$, and $D$ is the domain in $P$ bounded by $\Sigma$. In particular, when
$\Sigma$ is a round $(n-1)$-sphere of radius $\varrho$ it follows that
\beq
\label{est2}
0\leq|H_r|\leq\frac{1}{\varrho^r}.
\eeq
\end{theorem}
This estimate was first obtained in the case of constant mean curvature
($r=1$) by Barbosa in \cite{B1}.
\begin{proof}
Let $\xi$ be the unit vector normal to $P$. Then $\xi$ is a constant vector
field in \Rn, and therefore $Y=\xi$ is a Killing field in \Rn. On the other
hand we also have that $n_D=\pm\xi$, so that from \rf{estimativa} we
obtain
\[
n|H_r|\mathrm{vol}(D)=
\left|\oint_{\partial M}h_{r-1}\g{\xi}{\nu}^rds\right|\leq
\oint_{\partial M}|h_{r-1}|ds,
\]
which yields \rf{est1}.

In particular, when $\Sigma=\mathbb{S}^{n-1}(\varrho)$ is a round sphere
of radius $\varrho$, then we have that $\tau_i=-1/\varrho$ for every $i=1,
\ldots, n-1$, so that $h_{r-1}=(-1)^{r-1}/\varrho^{r-1}$. Besides, the
domain $D$ is an $n$-dimensional round ball of radius $\varrho$, with
volume
$n\ \mathrm{vol}(D)=\varrho\ \mathrm{vol}(\mathbb{S}^{n-1}(\varrho))$, and
the estimate \rf{est1} becomes \rf{est2}.
\end{proof}

Let us consider now the case of a hypersurface immersed into the hyperbolic space
\Hn. As in Section\rl{sectionflux}, it will be appropriate to use the Minkowski 
space model of \Hn, 
\[
\Hn=\{ x=(x_0,\ldots,x_{n+1})\in\mathbb{R}^{n+2}_1 : \g{x}{x}_1=-1, \quad x_0>0 \}.
\]
We may assume, up to an isometry of \Hn, that the totally geodesic 
hyperplane $P$ containing $\Sigma$ is given by
\[
P^n=\Hn\cap\{ x\in\mathbb{R}^{n+2}_1 : x_{n+1}=0 \}.
\]
In this case, the unit vector normal to $P$ in \Hn\ is given by
$\xi(p)=e_{n+1}=(0,\ldots,0,1)\in\mathbb{R}^{n+2}_1$ for every $p\in P$.
Observe that, for every arbitrary fixed point $a\in P$, the vector field given by
\[
Y(p)=-\g{p}{a}e_{n+1}+\g{p}{e_{n+1}}a, \quad p\in\Hn,
\]
is a Killing vector field on \Hn\ which is orthogonal to $P$, since at
every $p\in P$
\[
Y(p)=-\g{p}{a}e_{n+1}=\cosh{(\tilde{\varrho}(p))}\xi(p),
\]
where $\tilde{\varrho}(p)$ is the geodesic distance along $P$  between $a$ and 
$p$. Let $D$ be the compact domain $D$ bounded by $\Sigma$ in $P$, then $n_D=\pm\xi$ 
and from \rf{estimativa} we obtain
\beq
\label{esthip1}
n|H_r|\int_D\cosh{\tilde{\varrho}}\ dD=
\left|\oint_{\Sigma}h_{r-1}\cosh{\tilde{\varrho}}\g{\xi}{\nu}^rds\right|.
\eeq
Choose $a\in\mathrm{int}(D)$. Then $\min_{D}\cosh{\tilde{\varrho}}=\cosh{\tilde{\varrho}(a)}=1$, 
so that from \rf{esthip1} we conclude that
\begin{eqnarray}
\label{esthip2}
\nonumber n|H_r|\mathrm{vol}(D) & \leq & n|H_r|\int_D\cosh{\tilde{\varrho}}\ dD\leq 
\oint_{\partial M}|h_{r-1}|\cosh{\tilde{\varrho}}ds\\
{} & \leq & \max_{\Sigma}\cosh{\tilde{\varrho}}\oint_{\partial M}|h_{r-1}|ds.
\end{eqnarray}
In particular, when $\Sigma$ is a geodesic sphere in $P$ of geodesic 
radius $\varrho$ and $a$ is chosen to be the geodesic center of $\Sigma$, 
then $\tilde{\varrho}(p)=\varrho$ at every $p\in\Sigma$,  
$|h_{r-1}|=\coth^{r-1}(\varrho)$, and \rf{esthip2} simply becomes
\beq
\label{esthip3}
n|H_r|\mathrm{vol}(D)\leq
\cosh{\varrho}\coth^{r-1}(\varrho)\mathrm{vol}(\Sigma).
\eeq
Moreover, in this case $D$ is the geodesic ball in $P$ of radius $\varrho$ 
centered at $a$, that is, 
\[
D=\{ p\in P : 1\leq -\g{p}{a}<\cosh{\varrho} \},
\]
and
\[
\Sigma=\partial D=\{ p\in P : -\g{p}{a}=\cosh{\varrho} \}.
\]
Observe then that $\Sigma$ is in fact a round $(n-1)$-sphere of Euclidean radius
$\sinh\varrho$, and $D$ is a round $n$-dimensional ball of Euclidean radius $\sinh\varrho$.
Therefore, $\mathrm{vol}(\Sigma)=(n/\sinh\varrho)\mathrm{vol}(D)$ and \rf{esthip3}
simplifies to
\[
|H_r|\leq\coth^{r}(\varrho)
\]

We summarize this as follows.
\begin{theorem}
\label{estimate2}
Let $\Sigma$ be an orientable $(n-1)$-dimensional compact submanifold contained in a totally 
geodesic hyperplane $P\subset\Hn$, and let \xH\ be an orientable immersed compact connected
hypersurface with boundary $\Sigma=\psi(\partial\m)$ and constant $r$-mean curvature $H_r$,
$1\leq r\leq n$. Then
\[
0\leq |H_r|\leq\frac{C}{n\ \mathrm{vol}(D)}\oint_{\partial M}|h_{r-1}|ds.
\]
Here $h_{r-1}$ stands for the $(r-1)$-mean curvature of $\Sigma\subset P$, $D$ is the domain in 
$P$ bounded by $\Sigma$, and $C=\max_{\Sigma}\cosh{\tilde{\varrho}}\geq 1$, where 
$\tilde{\varrho}(p)$ is the geodesic distance along $P$ between a fixed arbitrary point 
$a\in\mathrm{int}(D)$ and $p$. 
In particular, when $\Sigma$ is a geodesic sphere in $P$ of geodesic radius $\varrho$, it 
follows that
\[
0\leq|H_r|\leq\coth^{r}\varrho.
\]
\end{theorem}

Finally, let us consider the case of a hypersurface immersed into the sphere
\Sn,
\[
\Sn=\{ x=(x_0,\ldots,x_{n+1})\in\mathbb{R}^{n+2} : \g{x}{x}=1 \}.
\]
We may assume, up to an isometry of \Sn, that the totally geodesic $n$-sphere $P$ containing 
$\Sigma$ is given by
\[
P^n=\Sn\cap\{ x\in\mathbb{R}^{n+2} : x_{n+1}=0 \}.
\]
In this case, the unit vector normal to $P$ in \Sn\ is given by
$\xi(p)=e_{n+1}=(0,\ldots,0,1)\in\mathbb{R}^{n+2}$ for every $p\in P$.
Observe that, for every arbitrary fixed point $a\in P$, the vector field given by
\[
Y(p)=\g{p}{a}e_{n+1}-\g{p}{e_{n+1}}a, \quad p\in\Sn,
\]
is a Killing vector field on \Sn\ which is orthogonal to $P$, since restricted to $p\in P$
\[
Y(p)=\g{p}{a}e_{n+1}=\cos{(\tilde{\varrho}(p))}\xi(p),
\]
where $\tilde{\varrho}(p)$ is the geodesic distance along $P$  between $a$ and 
$p$. Suppose that $\Sigma$ is contained in an open hemisphere $P_+$ of $P$ 
determined by an equator $S$ of $P$, and let $D$ be the compact domain $D$ bounded by $\Sigma$ 
in $P_+$. Then $n_D=\pm\xi$ and from \rf{estimativa} we 
obtain
\beq
\label{estesf1}
n|H_r|\left|\int_D\cos{\tilde{\varrho}}\ dD\right|=
\left|\oint_{\Sigma}h_{r-1}\cos{\tilde{\varrho}}\g{\xi}{\nu}^rds\right|.
\eeq
Choose $a\in\mathrm{int}(D)$. Since we are assuming that $\Sigma=\partial D$ is 
contained in the open hemisphere $P_+$, then $0\leq\tilde{\varrho}<\pi/2$ 
on $D$ and $\min_{D}\cos{\tilde{\varrho}}>0$, so 
that from \rf{estesf1} we conclude that
\begin{eqnarray}
\label{estesf2}
\nonumber n|H_r|\min_{D}\cos{\tilde{\varrho}}\mathrm{vol}(D) & \leq & 
n|H_r|\int_D\cos{\tilde{\varrho}}\ dD\leq 
\oint_{\partial M}|h_{r-1}|\cos{\tilde{\varrho}}ds\\
{} & \leq & \max_{\Sigma}\cos{\tilde{\varrho}}\oint_{\partial M}|h_{r-1}|ds.
\end{eqnarray}
In particular, when $\Sigma$ is a geodesic sphere in $P$ of geodesic 
radius $\varrho<\pi/2$ and $a$ is chosen to be the geodesic center of $\Sigma$, 
then $\tilde{\varrho}(p)=\varrho$ at every $p\in\Sigma$,  
$|h_{r-1}|=\cot^{r-1}(\varrho)$. Now a computation similar to that in 
hyperbolic space leads us from \rf{estesf2} to 
\[
|H_r|\leq\cot^{r}(\varrho),
\]
because, in this case $\Sigma$ is in fact a round $(n-1)$-sphere of Euclidean radius
$\sin\varrho$, and $D$ is a round $n$-dimensional ball of Euclidean radius $\sin\varrho$.
Summing up, we can state the following result.
\begin{theorem}
\label{estimate3}
Let $\Sigma$ be an orientable $(n-1)$-dimensional compact submanifold contained in an open 
totally geodesic hemisphere $P_+\subset\Sn$, and let \xS\ be an orientable immersed compact 
connected hypersurface with boundary $\Sigma=\psi(\partial\m)$ and constant $r$-mean curvature 
$H_r$, $1\leq r\leq n$. Then
\[
0\leq |H_r|\leq\frac{C}{n\ \mathrm{vol}(D)}\oint_{\partial M}|h_{r-1}|ds.
\]
Here $h_{r-1}$ stands for the $(r-1)$-mean curvature of $\Sigma\subset P$, $D$ is the domain in 
$P_+$ bounded by $\Sigma$, and $C=\max_{\Sigma}\cos{\tilde{\varrho}}/\min_D\cos{\tilde{\varrho}}$, 
where $\tilde{\varrho}(p)$ is the geodesic distance along $P_+$ between a fixed arbitrary point 
$a\in\mathrm{int}(D)$ and $p$. 
In particular, when $\Sigma$ is a geodesic sphere in $P_+$ of geodesic radius $\varrho<\pi/2$, 
it follows that
\[
0\leq|H_r|\leq\cot^{r}\varrho.
\]
\end{theorem}

\section{Symmetry for hypersurfaces in hyperbolic space}
\label{hyperbolic}
Hyperbolic space is rich in totally umbilic hypersurfaces. Besides the 
totally geodesic hyperplanes, there are the horospheres, the hyperspheres and 
the equidistant hypersurfaces. In all of them, the second fundamental form 
is proportional to the metric by a constant factor, and therefore they all 
have constant $r$-mean curvature, for $1\leq r\leq n$. After an 
appropriate choice of the unit normal vector field, hyperspheres have 
$r$-mean curvature bigger than 1, horospheres have $r$-mean curvature 1, 
and equidistant hypersurfaces have $r$-mean curvature in the interval 
$(0,1)$.

Let us fix a totally geodesic hyperplane $P^n\subset\Hn$ and a geodesic sphere 
$\Sigma^{n-1}\subset P^n$ in $\Hn$. Then each of the totally umbilic 
hypersurfaces above contains at least a compact domain $M^n$ with  
boundary being the sphere $\Sigma$. Those examples are called the 
\emph{spherical caps} in hyperbolic space. That terminology is due to the fact that, working 
in the half-space model of hyperbolic space, after an appropriate isometry of $\Hn$, 
the totally umbilic hypersurfaces above are given as intersections of \Hn\ 
with Euclidean spheres in \Rn. Because of the existence of these examples in 
\Hn, it is natural to consider the \emph{conjecture of the spherical cap} in hyperbolic space.

In this context, the corresponding result analogous to our Theorem\rl{calotaRn} for the case 
of hypersurfaces in hyperbolic space can be stated as follows.
\begin{theorem}
\label{simetriaHn}
Let $\Sigma^{n-1}$ be an strictly convex compact $(n-1)$-dimensional 
(connected) submanifold of a totally geodesic hyperplane $P^n\subset\Hn$, and let 
$M^n\subset\Hn$ be a compact (connected) embedded  hypersurface with boundary $\Sigma$. 
Let us assume that for a given $2\leq r\leq n$, the $r$-mean curvature 
$H_r$ of $M$ is a nonzero constant. Then $M$ has all the symmetries of $\Sigma$. 
In particular, when the boundary $\Sigma$ is a geodesic sphere in $P^n\subset\Hn$, then  
$M$ is a spherical cap.
\end{theorem}
As a consequence of Theorem\rl{simetriaHn} we can conclude, as in the 
Euclidean case, that the conjecture of the spherical cap is true for the 
case of embedded hypersurfaces with constant $r$-mean curvature in 
hyperbolic space, when $r\geq 2$.
\begin{corollary}
The only  compact embedded hypersurfaces in \Hn\ with constant $r$-mean curvature $H_r$
(with $2\leq r\leq n$) and spherical boundary are
\begin{itemize}
\item the geodesic balls of a totally geodesic hyperplane (with $H_r=0$);
\item the geodesic balls of an equidistant hypersurface (with $0<|H_r|<1$);
\item the geodesic balls of a horosphere (with $|H_r|=1$);
\item the geodesic balls of a hypersphere (with $|H_r|>1$).
\end{itemize}
\end{corollary}

\begin{proof}[Proof of Theorem\rl{simetriaHn}]
Let us work in the half-space model of hyperbolic space. We may assume, up 
to an isometry of \Hn, that the totally geodesic hyperplane $P$ is given 
by
\begin{equation}
P=\{x=(x_1,\ldots,x_{n+1})\in \mathbb H^{n+1};|x|=1,\; x_{n+1}>0\}.
\end{equation}
Let ${\mathcal B}$ be the connected component of $\Hn\backslash P$ containing the point 
$(0,\ldots,0,2)\in\Hn$. We will first see that there exists an interior 
elliptic point, that is, a point $p_0\in {\rm int}(M)$ where all the principal curvatures 
of \m\ are positive (after an appropriate orientation of \m). In fact, since $H_r$ is a nonzero constant, 
\m\ cannot be enterily 
contained in $P$. After an inversion with center $(0,\ldots,0)\in\Rn$ which fixes
$P$ (an isometry of \Hn), if necessary, we may assume that $M\cap\mathcal B\neq\emptyset$.
Let $C\subset P$ the geodesic sphere in $P$ given as the boundary of a geodesic 
ball in $P$ centered at the point $(0,\ldots,0,1)$ and containing 
$\Sigma$. Let us consider $\Gamma^{\varepsilon}\subset\Hn$ the equidistant 
sphere with center on the vertical geodesic through the center of $C$ such 
that $\Gamma^{\varepsilon}\cap P=C$, and such that the exterior angle 
between $\Gamma^{\varepsilon}$ and the asymptotic boundary of \Hn\ 
is $\frac{\pi}{2}-\varepsilon>0$. Since $\Gamma^{\varepsilon}\to P$ as $\varepsilon\to 0$, 
and taking into account that $\m\cap\mathcal B\neq\emptyset$, we may 
choose $\varepsilon>0$ such that $\Gamma^{\varepsilon}\cap M\neq\emptyset$. 
Besides, since $\Sigma$ is contained in the geodesic ball in $P$ bounded 
by $C$, then the points in $\Gamma^{\varepsilon}\cap M$ are interior points of \m.
Now, for every $t\geq 0$, let us consider $\Gamma^{\varepsilon}_t\subset\Hn$ the 
equidistant sphere in $\Hn$ obtained from $\Gamma^{\varepsilon}$ by an 
homothety centered at $(0,\ldots,0)\in\Rn$ (which is also an isometry 
of \Hn), and let us define $\Gamma^{\varepsilon}_0=\Gamma^{\varepsilon}$. If $t$ is large 
enough, then $\Gamma^{\varepsilon}_t$ englobes $M$; thus, we may find 
$t_0>0$ such that $M$ is tangent to $\Gamma^{\varepsilon}_{t_0}$ at a point $p_0$, 
which is necessarily an interior point of \m. Finally, it is easy to conclude 
that the normal curvatures of \m\ at $p_0$, with respect to the normal direction of the 
mean curvature vector of  $\Gamma^{\varepsilon}_{t_0}$, are greater or 
equal to those of $\Gamma^{\varepsilon}_{t_0}$, which are positive. In particular, 
choosing the appropriate orientation of \m, all the principal curvatures of \m\ at $p_0$ are 
positive.

Therefore, we may assume that $H_r=H_r(p_0)$ is a positive constant. This 
implies that for every $1\leq j\leq r-1$,  the Newton transformation $T_j$ 
is positive definite on \m\ (see \cite[Proposition 3.2]{BC}), and in 
particular the mean curvature is positive on \m, so that we may assume 
that \m\ is oriented by the mean curvature vector field. From Proposition \ref{transverse} we 
know that \m\ is transverse to $P$ along the boundary $\partial M$. This 
implies that, in a neighborhood of the boundary $\partial M$, \m\ is contained in 
one of the two connected components of $\Hn\backslash P$, which, without 
loss of generality, can be assumed to be $\mathcal B$. Beside, we may also 
assume that \m\ is globally transverse to $P$.

In this situation, we will prove that \m\ is above $P$, that is, $M\subset \bar {\mathcal B}$. 
Let us consider $\widetilde M$ the connected component of $M\cap \bar{\mathcal B}$ 
containing $\Sigma$. Then, $\widetilde M$ is a compact embedded hypersurface in \Hn\ with 
boundary $\partial\widetilde M$ contained in $P$. If the boundary $\partial\widetilde M$ were 
connected, then $\widetilde M=M$ and there is nothing to prove. Our 
objective is to show that actually $\partial\widetilde M$ must be 
connected. We will prove it by showing that assuming that 
$\partial\widetilde M$ is not connected yields a contradiction.

Thus, let us assume that the boundary $\partial\widetilde M$ consists of a finite 
number of disjoint connected compact emdedded $(n-1)$-dimensional submanifolds 
$\Sigma_i\subset P$ ($0\leq i\leq k$), with $\Sigma_0=\Sigma$. We orient 
this configuration as in Section \ref{configuration}, with $\widetilde M$ oriented 
by the mean curvature vector of \m. Let $\nu$ be the outward pointing 
conormal to $\widetilde M$ along each connected component of $\partial\widetilde M$. 
Then, the mean curvature vector of $M$, together with $\nu$, allows us to orient each $\Sigma_i$. 
Let $\eta$ be the unitary vector field normal to $\Sigma$ in $P$ which points outward with 
respect to the domain $D$ bounded by $\Sigma$ in $P$, and let $\xi$ be the
unique unitary vector field normal to $P$ which is compatible with $\eta$ and with the 
orientation of $\Sigma$. Now, there exists a unique choice for the unitary vector field 
$\eta_i$ normal to $\Sigma_i$ in $P$ which is compatible with the orientation of $\Sigma_i$ 
and with the orientation of $P$ given by $\xi$. We remark that we 
cannot ensure here that, for $i\geq 1$, $\eta_i$ points outward to the domain $D_i$ bounded
by $\Sigma_i$ in $P$. In this way, we have that formula (\ref{transverse}) holds at each point 
$p\in\partial\widetilde M$ with $r=1$, giving
\begin{equation}
\label{P_1nu}
\langle T_1\nu,\nu\rangle(p)=-s_1(p)\langle\xi,\nu\rangle(p)
\end{equation}
Here $s_1$ denotes the trace of the shape operator, with respect to $\eta_i$, of the inclusion 
$\Sigma_i\subset P$ which contains the point $p$.

As $\Sigma$ is a compact strictly convex submanifold of $P$ and $\eta$ points outward of $D$, 
then $s_1<0$ on $\Sigma$. On the other hand, as $T_1$ is positive definite on $M$, it follows 
from (\ref{P_1nu}) that 
$\langle\xi,\nu\rangle>0$ on $\partial M$.  Besides $\widetilde M\subset\mathcal{B}$ implies 
that $\langle\xi,\nu\rangle>0$ on each component of $\partial\widetilde M$. Hence, along 
$\Sigma$, the mean curvature vector of $M$ points into $D$. 
Therefore, if $\partial\widetilde M$ has a connected component contained 
in the interior of $D$, then there exists at least one component 
$\Sigma_i$, for some $i\geq 1$, contained in the interior of $D$ on which the mean curvature 
vector of $M$ points outward to the domain $D_i\subset P$ bounded by $\Sigma_i$ in $P$. As 
$\langle\xi,\nu\rangle>0$ on $\Sigma_i$, then $\eta_i$ must point into $D_i$. This contradicts 
the formula (\ref{P_1nu}), because if $\eta_i$ points into $D_i$, then we can easily conclude 
from the compactness of $\Sigma_i$ that there must be a point $p\in\Sigma_i$ where $s_1(p)>0$.
It then follows that the connected components of $\partial\widetilde M$ 
must be all contained in $P\backslash D$. 

Now, let us assume that there exists one of them, say $\Sigma_j$ ($j\geq 1$), 
which is homotopic to $\Sigma$ in $P\backslash D$. Without loss of generality, 
we may assume that, between $\Sigma_j$ and $\Sigma$  there is no other 
component of $\partial\widetilde M$ which is homotopic to $\Sigma$ in $P\backslash D$. 
We showed above that, along $\Sigma$, the mean curvature vector of $M$ points into $D$. 
Therefore, along $\Sigma_j$, the mean curvature vector of $M$ must point outward of
the domain $D_j\subset P$ bounded by $\Sigma_j$ in $P$. 
Since $\langle \xi,\nu\rangle>0$ on $\Sigma_j$, it then follows that the 
unitary vector field $\eta_j$ normal to $\Sigma_j$ in $P$ points into $D_j$. 
This situation gives again a contradiction with formula (\ref{P_1nu}), 
because if $\eta_j$ points into $D_j$, then there must be a point $p\in\Sigma_j$ where 
$s_1(p)>0$.

Finally, it only rests the case where $\partial\widetilde M$ has a 
connected component $\Sigma_l$ ($l\geq 1$) which is contained in 
$P\backslash D$ and is null homotpic in $P\backslash D$. However, this final
possibility is discarded by using the Alexandrov reflection 
technique \cite{Al}, exactly as in the proof of \cite[Theorem 1]{BEMR} or 
\cite[Theorem 3.1]{NR}. For the sake of completeness, we will include here 
the argument. Let $\gamma$ an infinite length geodesic in $P$ starting at a point 
of $D$ and intersecting $\Sigma_l$ in at least two points. 

Consider a family $Q(t)$, $t<\infty$, of geodesic hyperplanes of $\Hn$ 
orthogonal to $\gamma$, such that for each $q\in\gamma$ there exists 
exactly one $Q(t)$ which intersects $\gamma$ orthogonally at $q$. Each 
$Q(t)$ is orthogonal to $P$, so a hyperbolic simmetry through $Q(t)$ 
leaves $P$ and ${\mathcal B}$ invariant. Now we apply Alexandrov 
reflection method to $M$ (observe that this can be done because the 
equation $H_r=\mathrm{constant}>0$, under the existence of an elliptic 
point, is an elliptic equation \cite{Kor}). For $t$ large enough $Q(t)$
is disjoint from $M$. As $t$ decreases, there must exist a first point of 
contact of some $Q(t)$ with $M$. One continues to decrease $t$ and 
considers the symmetries of $M$ through the geodesic hyperplanes $Q(t)$. Sinces 
$\gamma$ intersects $\Sigma_l$ in at least two points, there must exist 
some hyperplane $Q(t_0)$ such that the symmetry of $M$ through $Q(t_0)$
will touch $M$ at an interior point. This occurs at an interior point 
since $\Sigma$ is convex and $\gamma$ intersects $\Sigma$ exactly at one point. 
Thus, $M$ is invariant under symmetry through $Q(t_0)$, which is impossible
(for $M$ would then be part of an embedded closed manifold with constant 
$r$-mean curvature, hence, a sphere. But a sphere cannot meet $P$ in more 
that one component).

Summing up, we conclude from the reasoning above that $\partial\widetilde M$ 
has no other connected component on $P$ except of $\Sigma$, and therefore 
$M\subset\bar{\mathcal B}$. Now that we know that $M$ is above $P$ and 
transverse to $P$ along $\partial M$, the proof finishes applying again 
the Alexandrov reflection method to $M\cup D$, exactly at in the final 
step of the proof of \cite[Theorem 2.1]{NR}.
\end{proof}

\section{Symmetry for hypersurfaces in sphere}
\label{sphere}
The totally umbilic hypersurfaces of  \Sn\ are given by the 
intersections of \Sn\ with the hyperplanes of Euclidean space $\mathbb R^{n+2}$. 
When the hyperplane passes through the origin of $\mathbb R^{n+2}$, they 
are totally geodesic, and when the hyperplane is an affine hyperplane, they are 
totally umbilic. We will refer to them as totally geodesic $n$-spheres and 
totally umbilic $n$-spheres of \Sn, respectively. They all have constant 
$r$-mean curvature. After an appropriate choice of the unit normal vector 
field, the totally umbilic $n$-spheres have $r$-mean curvature 
$H_r=\cot^r(\varrho)$, where $\varrho>0$ denotes the geodesis radius of the convex geodesic 
ball of $\Sn$ whose boundary is the totally umbilic $n$-sphere.

Let us fix a totally geodesic $n$-sphere $P^n\subset\Sn$ and a geodesic sphere 
$\Sigma^{n-1}\subset P^n$ in $\Sn$. Then, fixed a value for $H_r$, there are two compact domains 
$M_1^n$ and $M_2^n$ of a totally umbilic $n$-sphere of \Sn\ whose boundaries are the geodesic 
sphere $\Sigma$. 
These examples are called the \emph{spherical caps} in \Sn. As in hyperbolic space, because of 
the existence of these examples in \Sn, it is also natural to consider the 
\emph{conjecture of the spherical cap} in \Sn. In this context, the 
corresponding result for the case of hypersurfaces in \Sn\ can be stated 
as follows.
\begin{theorem}
\label{simetriaSn}
Let $\Sigma^{n-1}$ be a convex $(n-1)$-dimensional submanifold of a totally geodesic $n$-sphere 
$P^n\subset\Sn$, and let $M^n\subset\Sn$ be a compact (connected) embedded hypersurface with 
boundary $\Sigma$. Let us assume that $M$ is contained in an open hemisphere \Snp,
and that the $r$-mean curvature $H_r$ 
of $M$ is a nonzero constant, for a given $2\leq r\leq n$.  
Suppose that the convex disc $D$ bounded by $\Sigma$ in $P$ contains a 
focal point of $P_1\cap P$, where $P_1=\partial\Snp$.
Then $M$ has all the symmetries of 
$\Sigma$. In particular, when the boundary $\Sigma$ is a geodesic sphere in $P^n\subset\Sn$, 
then $M$ is a spherical cap.
\end{theorem}
\begin{corollary}
Let \m\ be a compact (connected) embedded hypersurface in \Snp\ with constant $r$-mean curvature 
$H_r\neq 0$ (with $2\leq r\leq n$) and spherical boundary contained in a 
totally geodesic $n$-sphere $P^n\subset\Sn$. Suppose that the convex disc $D$ bounded by 
the spherical boundary of \m\ in $P$ contains a focal point of $P_1\cap 
P$, where $P_1=\partial\Snp$. Then \m\ is a spherical cap.
\end{corollary}
As we had pointed out before, corresponding results for $r=1$ can be found in \cite{Lira}.

Before go further, it is needed to fix a suitable notion of symmetry in the spherical space form. 
This is done in the definition below.

\begin{definition} We say that a totally geodesic $n$-sphere $Q$ is a $n$-\emph{sphere of 
symmetry} of a subset $S$ of $\mathbb{S}^{n+1}$ if for each point $p\in
S$ and any complete 
geodesic $ \gamma$  perpendicular to $Q$ and containing $p$, we have $\tilde p\in S$, where
$\tilde{p}$ is the
point of $\gamma$ such that $p$ and $\tilde{p}$ lie in opposite hemispheres of
$Q$ at distance less than or equal to $\frac{\pi}{2}$ and $\mathop{\rm{dist}}(\tilde
{p},Q)=\mathop{\rm{dist}}(p,Q)$.
\end{definition}

We observe that the choice of $D$ in this section is compatible with the orientations established in the Section 4, what allows us to use the  calculations made in the earlier parts of the article.

\begin{proof}

Let $a\in\Sn$ and consider $P_1=\{x\in\Sn :\langle x,a\rangle =0\}$ the totally geodesic 
$n$-sphere
which defines the open hemisphere $\Snp=\{x\in\Sn :\langle x,a\rangle >0\}$ where $M$ is 
contained. Now, we may assume without loss of generality that the totally 
geodesic $n$-sphere containing the boundary of $M$ is 
\[
P=\{x\in\Sn:\langle x,e_0\rangle=0\}, a\neq e_0.
\]
Our first objective is to see that there exists an interior elliptic point of $M$, that is, a 
point $p_0\in {\rm int}(M)$ where all the 
principal curvatures of $M$ have the same sign. To see it, let $B_t(a)\subset\Snp$ be the 
geodesic 
ball with center $a$ and geodesic radius $t$, where $0<t<\pi/2$, and let 
$S_t(a)=\partial B_t(a)$ be the corresponding geodesic sphere. Since $M$ 
is compact and $M\subset\Snp$, there exists a minimum value $t^\prime$ such that 
$M\subset\overline{B_{t^\prime}(a)}$, and a contact point $p_0\in M\cap S_t(a)$. 
Observe that the height function $\langle x,a\rangle$ on $M$ attains its minimum value precisely 
at that contact point. Therefore, if such a contact point is an interior point of 
$M$, then it is also a tangency point and all the principal curvatures of 
$M$, with respect to the unit normal vector field of $S_{t^\prime}(a)$, 
are positive at $p_0$. If the contact point is a boundary point, then we 
can consider a geodesic ball $B_t(a)$ with $t>t'$ so that 
$B_t(a)\cap\Sigma=\emptyset$. Now, we can simultaneously move the center $a$ of the geodesic 
ball and decrease its radius, keeping always $M$ contained in the interior 
of this geodesic ball, and we consider the intersection of this geodesic 
ball with \Snp. From this process it follows that either some geodesic 
ball $B_t(a')\cap\Snp$ is tangent to $M$ at an interior point, or $M$ is 
entirely contained in the totally geodesic $n$-sphere $P$. However, the 
second possibility cannot happen because $H_r$ is a nonzero constant. Then, reasoning as 
above, such an interior tangency point is an elliptic point of $M$. Thus, we may always 
(including when $r=$even) assume that the $r$-mean curvature $H_r$ of $M$ is a positive constant. 
This implies that $T_j$ is positive definite on $M$, for each $1\leq j\leq r-1$, and, since 
$H=H_1>0$, we may orient $M$ by the mean curvature vector. By the Proposition \ref{transverse}, 
we conclude that $M$ is transversal to $P$ along its boundary $\partial M$. So, there is a 
neighbourhood $\mathcal{U}$ of $\Sigma$ in $M$ contained in only one of the hemispheres  
$\overline{P^+}$ and $\overline{P^-}$ determined by $P$. We fix 
$\mathcal{U}\subset\overline{P^+}$. 

Let $D\subset P$ be the domain bounded by $\Sigma$ which does not contain points of 
$\Sigma_1:=P_1\cap P$. Denote by ${\rm int}(D)$ the interior of $D$ in $P$ and by 
${\rm ext}(D)$ the subset $P-D$. According to this notation, we have

\begin{claim}
If $M\cap{\rm int}(D)\neq\emptyset$, then $M\cap{\rm ext}(D)\neq\emptyset$. 
\end{claim}

We omit the proof of this claim since it follows the same guidelines from the similar one for 
the hyperbolic case (see Section 10).

For guarantee that $M\cap P=\Sigma$ it suffices then, by the Claim, to prove that 
$M\cap {\rm ext}(D)=\emptyset$. Suppose otherwise, that is, suppose that 
$M\cap\textrm{ext}(D)\neq\emptyset$. So, we may assume, without loss of generality, that 
$M\cap\textrm{ext}(D)$ consists of a finite
number of disjoint connected embedded submanifolds of $P$.

In order to apply correctly the reflection procedure,  we consider now on $M$ the connected 
component containing $\Sigma$ of the topological embedded submanifold $N$ obtained after 
excision of small annuli in $M$ surrounding each one of the domains $D_i$ bounded by the 
components of $\Sigma\cap {\rm int} D$ containing no points of $\Sigma$  and gluing domains 
homeomorphic to each $D_i$ at the boundary of these annuli (see \cite{BEMR} for the similar 
device in $\mathbb{R}^3$). This construction allows us to consider $M\cup D$ separating 
$\mathbb{S}^{n+1}$ in two connected components. By $\Omega$ we denote the component that 
contains no points of $\Sigma_1$. Note that the set $M\cup \mathop{\rm{ext(D)}}$ is not 
diminished in this process. In fact, the only components of $M\cup \mathop{\rm{ext(D)}}$ that 
could be discarded are the ones that contains points in $\mathop{\rm{int(D)}}$ and points in 
$\mathop{\rm{ext(D)}}$, whose existence should  oblige the mean curvature vector  to point 
outside $\Omega$, contradicting the Maximum Principle applied to geodesic graphs in $\Sn$ 
(see \cite{JRF}).  

\noindent{\bf Case 1.} Suppose initially that there exist components $\Sigma_{k}$ of 
$M\cap\textrm{ext}(D)$ homologous to zero in
$P-\textrm{int}(D)$. For each $k$, denote by $M_k$ the connected component of $M$ which has 
boundary $\Sigma_{k}$ and contains no points of $\Sigma$. We note that $M_k$ contains points of 
$\overline{P^-}$ in a neighbourhood of $\Sigma_k$.

We fix $\Sigma_1=\{x=(0,0,x_2,\ldots,x_{n+1})\in\mathbb{S}^{n+1}\}$ and $x_1>0$ throughout
the hemisphere of $P$ which contains no points of $\Sigma$. Define $P_1(t),\,
0\le{t}\le\pi$, as the family of totally geodesic $n$-spheres such that
$P_1(t)\cap P=\Sigma_1$, for all $t$, and $P_1(\alpha)=P_1$, where $\alpha$ is the angle
between $P_1$ and $P$. The normal vector to $P_1(t)$ is given by
$n_t=(\cos t,-\sin t,\ldots,0,0)$.

For $t,\,0\le{t}<{\alpha}$, let $M_k^{-}(t)$ be the set
$\{x\in{M};{\langle{x,n_{\alpha}}\rangle>0}\,{\mathrm{and}}\,{\langle{x,n_{t}}\rangle<0}\}$
and let \ $\widetilde{M}_k(t)$ be  the reflected image of \ $M_k^{-}(t)$ \ through \
$P_1(t)$, \ i.e., \
$\,\widetilde{M}_k(t)=\{\tilde{x}\in{\mathbb{S}^{n+1}};\tilde{x}={x-2\langle{x,n_{t}}\rangle{n_t}},\,x\in{M_k^{-}(t)}\}$.

By the fact that ${P_1(\alpha)}\cap{M_k}={\emptyset}$, there exists $t_0,\,0\le{t_0}<{\alpha}$, 
such that
\begin{itemize}
\item[(i)] ${P_1(t_0)}\cap{M_k}\neq{\emptyset}$;
\item[(ii)] ${P_1(t)}\cap{M_k}={\emptyset}$, for all $t>t_0$.
\end{itemize}

So, $M_k$ is tangent to $P_1(t_0)$ at their common points and there is a neighbourhood of
each one of these points in $M_k$
which is a geodesic graph over a domain in $P_1(t_0)$.

Thus, unless ${t_0}=0$, we have that
$\widetilde{M}_k(t)\subset\overline{P^-}$ for $t$
sufficiently close to
$t_0$.
However, for ${t_0}= 0$, we consider a rotation of $P$ by a small angle, fixing ${\Sigma}_1$, 
to return to the previous situation.

We claim that $M_k$ is a geodesic graph over the domain $D_k$ in $P$ bounded by $\Sigma_k$  
which contains no points of $\Sigma_1$, with $\widetilde{M}_k(0)\subset\textrm{int}(\Omega)$. 
In particular, $M_k$ is not perpendicular to $P$ at points of $\Sigma_k$. Otherwise, there 
exist $k$ and $t_1\in [0,t_0)$ for which occurs one of the following possibilities:
\begin{itemize}
\item[(i)] ${\widetilde{M}_k(t_1)}\cap{M}$ contains interior points of $M_k$;
\item[(ii)] $P_1(t_1)$ is perpendicular to $M_k$ at points of $\partial M_k(t_1)$;
\item[(iii)] ${\widetilde{M}_k(t_1)}\cap{M}$ contains points of
$(M-\Sigma)-M_k$;
\item[(iv)] $\widetilde{M}_k(t_1)$ contains points of $\Sigma$.
\end{itemize}

The cases (i) to (iii) are all ruled out by the Maximum Principle. In fact, otherwise,
$P_1(t_1)$ should be a sphere of symmetry and $M$ a compact hypersurface
without boundary (see \cite{Ko}, p. 572-573).

Thus, we conclude that the points of $\widetilde{M}_k(t_1)$ away from $\Sigma$ are
contained in $\textrm{int}(\Omega)$.
Since there exists a neighbourhood $\mathcal{U}$ of $\Sigma$ in $M$ as above, the reflected 
image of
$\widetilde{M}_k(t_1)$ through
$P=P(\pi)$ is not contained in $\Omega$, if we suppose $\widetilde{M}_k(t_1)\cap\Sigma\neq\emptyset$. In particular, the reflection of $\widetilde{M}_k(t_1)$ through $P$ is not contained in the open domain bounded by $\widetilde{M}_k(t_1)$ in $\textrm{int}(\Omega)$. Therefore, a sphere
$P_1(\tau),\,\alpha<\tau<\pi$, should exist such that the reflected image
of $\widetilde{M}_k(t_1)$ through $P_1(\tau)$ is tangent to $\widetilde{M}_k(t_1)$ and, in this 
way, $P_1(\tau)$ is sphere of symmetry of $\widetilde{M}_k(t_1)$. Hence, since the portion of 
$\widetilde{M}_k(t_1)$ lying between $P_1(\alpha)$ and $P_1(\tau)$ does not contain points of 
$\partial \widetilde{M}_k(t_1)=\partial M_k(t_1)$ in $P_1(t_1)$, we obtain a contradction, 
proving that the case (iv) does not occur at any instant $t\in [0,t_0)$.

Notice that it is equally
impossible to have $\widetilde{M}_k(t_1)$ tangent to $M$ at points of opposite
orientation, because if
it is the case, then the reflected image of a portion of $M_k$ would have left 
$\textrm{int}(\Omega)$ before $t_1$.

So, we conclude from the impossibility of the cases (i) to (iv), for each $t\in [0,t_0)$, that 
$\overline{M_k(t)}$ is a geodesic graph over the domain in $P_1(t)$ bounded by 
$\partial M^-_k(t)$ which does not contain points of $\Sigma_1$, with 
$\widetilde{M}_k(t)\subset\textrm{int}(\Omega)$, proving the claim. Besides this, we guarantee 
that $\overline{M_k^-(t)}$ is not perpendicular to $P_1(t)$ at any point of $\partial M_k^-(t)$.

\noindent {\bf Case 2.} Now, suppose there exist components of $M\cap \textrm{ext}(D)$ 
homologous to $\Sigma$. We will prove that whenever exist such components, they are graphs over 
domains in $P$.

This case is handled as in \cite{Lira}, with minor modifications concerning the utilization of 
the flux formula there, which must be changed by the appropriate formula \ref{fluxo3}.

Now, as above, define for $t\in (\alpha,\frac{\pi}{2}]$ the submanifold of $M$ given by 
$M^-(t)=\{x\in P^+;\langle x,n_\alpha\rangle>0\,\mathrm{and},\langle x,n_t\rangle<0\}$
and its reflected image through $P_1(t)$ as
$\widetilde{M}(t)=\{\tilde{x}\in{\mathbb{S}^{n+1}};\tilde{x}={x-2\langle{x,n_{t}}\rangle{n_t}},\,x\in{M^{-}(t)}\}$.

Since $M\cap{P_1}=\emptyset$, either $M$ is contained in the open hemisphere
determined by $P_1(\frac{\pi}{2})$ containing $\Sigma$, or exists 
$t_0\in (\alpha,\frac{\pi}{2}]$ such that
\begin{itemize}
\item[(i)] $P_1(t_0)\cap M\neq\emptyset$;
\item[(ii)] $P_1(t)\cap M=\emptyset$, for all $\alpha<t<t_0$.
\end{itemize}

For $t_0=\frac{\pi}{2}$, there is a neighbourhood of $M$ that is a graph over a domain of 
$P_1(\frac{\pi}{2})$.

If $t_0<\frac{\pi}{2}$, suppose that there exists $t_1\in (t_0,\frac{\pi}{2}]$, for which holds 
one of the statements below:
\begin{itemize}
\item[(i)] $\widetilde{M}(t_1)$ is tangent to $M$ at interior points;
\item[(ii)] $P_1 (t_1)$ is perpendicular to $M$ at some points of  $M\cap P_1(t_1)$;
\item[(iii)] $\widetilde{M}(t_1)\cap\Sigma\neq\emptyset$.
\end{itemize}

If (i) or (ii) occurs, then $P_1(t_1)$ is a sphere of symmetry of $M$. However,
$\Sigma$ is contained in only
one of the hemispheres determined by $\Sigma_1=P_1(t_1)\cap P$ on $P$.

Suppose (iii) occurs; if exists $p\in M^-(t_1)$ such that $\tilde{p}\in\Sigma$, then $p$ and 
$\tilde{p}$
are points at the same distance from $P_1(t_1)$ on a geodesic $\Sigma$ perpendicular to 
$P_1(t_1)$. If
$t_1=\frac{\pi}{2}$, then $p\in{P}$, since $P$ is totally geodesic. If $t_1<\frac{\pi}{2}$, we 
have $\textrm{dist}(\tilde{p},P_1(t_1))<\frac{\pi}{2}$; thus, 
$\textrm{dist}(p,\tilde{p})<2t_1<\pi$ what implies $p\in P^-$. Both
situations contradict the fact that $M^-(t_1)\subset P^+$.

We conclude that $\widetilde{M}(t_1)\subset\textrm{int}(\Omega)$,
for all $t\in(\alpha,\frac{\pi}{2}]$.
Furthermore, $\overline{M^-(\frac{\pi}{2})}$ is a geodesic graph over the domain bounded by 
$\partial M^-(\frac{\pi}{2})$ in $P_1(\frac{\pi}{2})$ containing points of $\Omega$.

Let $p\in D$ the geodesic center of $\Sigma_1$ and $\sigma$ an arc of geodesic starting from $p$ 
passing through $\Sigma$ and crossing orthogonally $\Sigma_1$. We may assume,
initially, that the component $\Sigma_{k_0}=\partial M_{k_0}$ of $M\cap\textrm{ext}(D)$
nearest from $\Sigma_1$ in the direction given by $\sigma$ is homologous 
to zero. Modifying slightly the direction of $\sigma$, if necessary, we may
assume that $\sigma$ crosses $\Sigma_{k_0}$ at least twice.

For each point $\sigma(t)$, $0\le t\le d$, we consider the intersection 
$Q(t)$ of $\mathbb{S}^{n+1}$ and the
Euclidean hyperplane containing the origin of $\mathbb{R}^{n+2}$ and perpendicular to
$\{x_0=0\}$ whose normal vector is $\sigma'(t)$. Denote by $Q^-(t)$ the hemisphere determined by 
$Q(t)$ containing $\sigma[t,d]$ and by $\mathcal{Q}_t$ the reflection through $Q(t)$. 

As we have proved, the portion of $M$ in the hemisphere $Q^-(d)$ determined by
$Q(d)=P_1(\frac{\pi}{2})$, if it is not empty,
is a geodesic graph over a domain in $P_1(\frac{\pi}{2})\cap P^+$ at distance less than 
$\frac{\pi}{2}$ from the sphere $P_1(\frac{\pi}{2})$. This remains true,
for $t$ sufficiently close to $d$. By the choice of $\Sigma_{k_0}$, we have that the first point 
of contact, if exists, between the planes $Q(t)$ and $M\cap\overline{P^-}$ must be in
$\Sigma_{k_0}$. More precisely, there exists $t_0\in (0,d)$ such that we have:
\begin{itemize}
\item[(i)] $M\cap Q^-(t)$ is contained in the portion of the geodesic cylinder over a domain of 
$Q(t)$ contained in $P^+$ and $\mathcal{Q}_t(M\cap Q^-(t))\subset\Omega\cap P^+$, for all 
$t<t_0$; 
\item[(ii)] $M\cap\overline{P^-}\cap Q(t_0)$ is a non-empty subset of $\Sigma_{k_0}$.
\end{itemize}

These statements follow from the fact that $M_{k_0}$ is, as proved above, contained in the 
geodesic cylinder over a domain in $P$ and $Q(t_0)$ is a totally geodesic sphere perpendicular 
to $P$. So, since $\sigma$ crosses $\Sigma_{k_0}$ at least twice, and $M_{k_0}$ is compact, 
there exists $t_1\in (0,t_0]$ so that $M\cap Q^-(t)$ is contained in the geodesic cylinder over 
a domain of $Q(t)$, in such a way that $\mathcal{Q}_t(M\cap Q^-(t))\subset\textrm{int}(\Omega)$, 
whenever $t>t_1$. Furthermore, one of the following assertions holds: 

\begin{itemize}
\item[(i)] $\mathcal{Q}_{t_1}(M\cap Q^-(t_1))$ is tangent to $M_{k_0}$ at
points not
belonging
to $Q(t_1)$ with the same orientation;
\item[(ii)] $Q(t_1)$ is perpendicular to $\overline{M\cap Q^-(t_1)}$ at
points of $Q(t_1)$ or,
equivalently, $\overline{\mathcal{Q}_{t_1}(M\cap Q^-(t_1))}$ is tangent
to $M$ at points of
$Q(t_1)$.
\end{itemize}

In any case, $Q(t_1)$ should be a sphere of symmetry of $M$. Let $p'$ be the last point of
$\Sigma$ in $\sigma[0,d)$. The distance between $p'$ and $Q(t_1)$ is less than $t_1$.
Thus, prolonging $\sigma$ until the point $\mathcal{Q}_{t_1}(p')$, we obtain an arc of geodesic 
of lenght strictly less than $2t_1<\pi$. Since $Q(t_1)$ is a sphere of
symmetry of $M$ and, in particular, of $\Sigma$, 
we have that $\mathcal{Q}(p')$ is a point of $\Sigma$. However, since that $\Sigma$ is convex, 
$\sigma$ does not return to $\Sigma$ until it has just crossed all of the hemisphere determined 
by $\Sigma_1$ in $P$ which does not contain points of $\Sigma$, that is, just after $t>\pi$. As 
$2t_1<\pi$, we have a contradiction. From this contradiction, we conclude that there is no 
components of $M\cap\textrm{ext}(D)$ homologous to zero outside the region in $P-\Sigma_1$ 
bounded by $\Sigma$ and some component of $M\cap\textrm{ext}(D)$ homologous to $\Sigma$; 
otherwise, there exists at least a direction $\sigma$ starting from $p$ so that the component 
of $M\cap\textrm{ext}(D)$ nearest from $\Sigma_1$ in its direction is homologous to a constant.

Now, suppose $\Sigma_{k_0}$ is homologous to $\Sigma$. By construction, it is clear that 
$\sigma$ crosses each component $\Sigma_k$ of $M\cap\textrm{ext}(D)$ homologous to $\Sigma$ at 
least once. So, proceeding as in \cite{Lira} we find a sphere of symmetry of $M$ before  reach 
$\Sigma$, a contradiction.

At the moment, we have proved that $M$ is contained in $\overline{P^{+}}$ and that
${M\cap{P}}=\Sigma$. Furthermore, $M^{-}(\frac{\pi}{2})$, if it is not empty, is a geodesic 
graph over a domain in $P_1(\frac{\pi}{2})$ having height less than $\frac{\pi}{2}$ and
its reflected image through $P_1(\frac{\pi}{2})$ is entirely contained in
$\textrm{int}(\Omega)$. 

Let then $R$ be a sphere of symmetry of $\Sigma$ and $q\in R\cap D$.
Let $\mu$ be the geodesic perpendicular to $R$ starting from $q$ and reaching a point 
$q'\in\Sigma_1$. We define $R(t)$, $0\le t\le 2\pi$, as the intersection of $\mathbb{S}^{n+1}$ 
and the Euclidean hyperplane containing the origin of $\mathbb{R}^{n+2}$ and perpendicular to 
$\{x_0=0\}$, whose normal is $\mu'(t)_{(0,\ldots,0)}$. It is clear that $R(0)=R$. Suppose that 
$R(d)$ and $P_1(\frac{\pi}{2})$ coincide. Then, the facts above imply that we have no touching 
points until the time $t=d$ on the reflection process through the spheres $R(t)$. However, 
since $M$ is compact and ${M\cap{P}}=\Sigma$, there exists $t_1\in{\lbrack{0,d})}$ such that 
$R(t_1)$ is a sphere of symmetry of $M$ and, in particular, of $\Sigma$. Since $R$ and $R(t_1)$ 
are both perpendicular to $\mu$, it follows from the convexity of $\Sigma$ that $R={R(t_1)}$, 
i.e., that $R$ is a sphere of symmetry of $M$.

If $R(d)$ and $P_1(\frac{\pi}{2})$ \ are \ distinct \ spheres, \ let
$\Sigma_2=P_1(\frac{\pi}{2})\cap R(d)$ and consider the totally geodesic spheres 
$T(t),\,0\le t\le\alpha_0$, obtained by rotation, fixing $\Sigma_2$, of $P_1(\frac\pi 2)$ 
towards $R(d)$, with $T(0)=P(\frac\pi 2)$ and $T(\alpha_0)=R(d)$. It is clear that
$\Sigma_2=\cap_t T(t)$. Moreover, we have that
$T(t)\cap\Sigma=\emptyset$, for all $t$, since each $T(t)$ is contained in the domain 
$\mathcal{C}$ of $\mathbb{S}^{n+1}$ bounded
by $P_1(\frac{\pi}{2})$ and $R(d)$ that does not contain points of $\Sigma$. Denote $T^-(t)$ and
$\mathcal T(t)$ as before. 

By continuity, we have that, for $t$ close enough to $0$, each component of $M\cap T^-(t)$, 
when this set is not empty, is  still a geodesic graph over a domain $T(t)$ at distance from 
$T(t)$ strictly less than $\frac{\pi}{2}$. Furthermore, since $M$ is compact, it is
possible to consider $t$ sufficiently small so that
$\mathcal T(t)(M\cap T^-(t))\subset\textrm{int}(\Omega)$. Thus, either this remains
true for each
$t\in (0,\alpha_0\rbrack$, or there exists $t_1\in (0,\alpha_0\rbrack$ such that one of the
following situations occurs:
\begin{itemize}
\item[(i)] $\mathcal T_{t_1}(M\cap T^-(t_1))$ is tangent to $M$ at points not
belonging
to $T(t_1)$ with the same orientation;
\item[(ii)] $T(t_1)$ is perpendicular to $\overline{M\cap T^-(t_1)}$ at points of
$T(t_1)$ or, equivalently, $\overline{\mathcal{T}_{t_1}(M\cap T^-(t_1))}$ is tangent to $M$ at
points of
$T(t_1)$.
\end{itemize}

In these cases, $T(t_1)$ should be a sphere of symmetry of $M$ and, in
particular, of $\Sigma$.
However, this contradicts the fact that there are no points of $\Sigma$ in $\mathcal{C}$. 
Therefore, we conclude from this contradiction that 
${{\mathcal{T}}_{\alpha_0}}(M\cap{T^{-}(\alpha_0)})$ is contained in $\Omega$ and that 
$M\cap{T^{-}}(\alpha_0)$ is either empty or a graph over ${T(\alpha_0)}={R(d)}$. In
this way, we return to the previous case. The theorem is proved. 
\end{proof}

\section*{Acknowledgements}
This paper is part of the PhD thesis of the third author \cite{Ma}, which 
was defended in Universidade Federal do Cear\'a, Fortaleza, Brazil, in February 2001. He 
would like to thank his adviser, Prof. Luqu\'esio P.M. Jorge, for his guidance.

This work was started while the first author was visiting the IHES at Bures-sur-Yvette, France, 
in Summer 2000, and it was also developed while he was visiting again the IHES in Summer 2002. 
This work also benefited from two visits of the first author to the Departamento de 
Matem\'atica of the Universidade Federal do Cear\'a, Fortaleza, Brazil, in November 2000 and 
April 2002. He would like to thank both institutions for their hospitality.

\end{document}